\documentclass[10pt]{amsart}

 \usepackage{graphicx}
 \usepackage[all]{xy}
 \usepackage{mathbbol}
 \usepackage{amscd,amsfonts,amssymb,amsmath,amsthm,latexsym}
\usepackage{mathtools}
 \usepackage{tikz-cd}
 \usepackage{color}
 \usepackage{hyperref}
\usepackage{extarrows} 
\usepackage{verbatim}

\makeatletter
\newcommand{\doublewidetilde}[1]{{%
  \mathpalette\double@widetilde{#1}%
}}
\newcommand{\double@widetilde}[2]{%
  \sbox\z@{$\m@th#1\widetilde{#2}$}%
  \ht\z@=.9\ht\z@
  \widetilde{\box\z@}%
}
\makeatother

\begin{document}

\newcommand{\suchthat}{\mid }
\newcommand{\bbC}{\mathbb{C}}
\newcommand{\bbZ}{\mathbb{Z}}
\newcommand{\field}{\mathbb{k}}
\newcommand{\pathalg}[1]{\field\langle #1\rangle}
\newcommand{\compalg}[1]{\field\langle\hspace{-0.075cm}\langle #1\rangle\hspace{-0.075cm}\rangle}
\newcommand{\idealM}{\mathfrak{m}}
\newcommand{\jacobalg}[1]{\mathcal{P}(#1)}
\newcommand{\coker}{\operatorname{coker}}
\newcommand{\image}{\operatorname{im}}
\newcommand{\End}{\operatorname{End}}
\newcommand{\Hom}{\operatorname{Hom}}
\newcommand{\Gr}{\operatorname{Gr}}
\newcommand{\Mod}{\mathbf{Mod}}
\newcommand{\DecMod}{\mathbf{DecMod}}
\newcommand{\soc}{\operatorname{soc}}
\newcommand{\rad}{\operatorname{rad}}
\newcommand{\myid}{1\hspace{-0.075cm}1}
\newcommand{\surf}{(\varsigma,\mathbb{M})}
\newcommand{\charac}{\varrho}
\newcommand{\sinis}{\operatorname{left}}
\newcommand{\destra}{\operatorname{right}}
\newcommand{\Tr}{\operatorname{Tr}}
\newcommand{\ZZ}{\mathbb{Z}}
\newcommand{\QQ}{\mathbb{Q}}
\newcommand{\calA}{\mathcal{A}}
\newcommand{\calG}{\mathcal{G}}
\newcommand{\calT}{\mathcal{T}}
\newcommand{\Match}{\mathrm{Match}}

\newcommand{\calB}{\mathcal{B}}

\numberwithin{equation}{subsection}

\theoremstyle{plain}
    \newtheorem{theorem}{Theorem}[section]
    \newtheorem*{teo*}{Theorem}
    \newtheorem{definition}[theorem]{Definition}
    \newtheorem{lemma}[theorem]{Lemma}
    \newtheorem{proposition}[theorem]{Proposition}
    \newtheorem{corollary}[theorem]{Corollary}
    \newtheorem{question}[theorem]{Question}
       \newtheorem{notation}[theorem]{Notation}

\theoremstyle{remark}
    \newtheorem{remark}[theorem]{Remark}
    \newtheorem{example}[theorem]{Example}
\newtheorem{case}{Case}
\newtheorem{subcase}{Subcase}[case]

\title[Skew-symmetrizable cluster algebras from surfaces and symmetric quivers]{Skew-symmetrizable cluster algebras from surfaces and symmetric quivers}
\author[Azzurra Ciliberti]{Azzurra Ciliberti}
\address{Fakult\"{a}t f\"{u}r Mathematik, Universit\"{a}tsstraße 150, D-44780 Bochum, Ruhr-Universit\"{a}t Bochum}\email{\href{mailto:azzurra.ciliberti@ruhr-uni-bochum.de}{azzurra.ciliberti@ruhr-uni-bochum.de}}

\begin{abstract}
\noindent 
We study skew-symmetrizable cluster algebras $\mathcal{A}$ associated with unpunctured surfaces $\Tilde{\mathbf{S}}$ endowed with an orientation-preserving involution $\varsigma$. We give a geometric realization of such cluster algebras by showing that cluster variables of $\mathcal{A}$ correspond to non-crossing $\varsigma$-orbits $[\gamma]$ of arcs of $\Tilde{\mathbf{S}}$, while clusters are given by admissible $\varsigma$-invariant triangulations. We establish a formula expressing some cluster variables of $\mathcal{A}$ in terms of those of a skew-symmetric cluster algebra of the same rank, which is combinatorially derived from $\mathcal{A}$.
We use this result to provide a cluster expansion formula in terms of perfect matchings of some labeled modified snake graphs constructed from the arcs of $[\gamma]$. Then, we associate a symmetric finite-dimensional algebra $A$ to $\mathcal{A}$, such that non-initial cluster variables correspond to some orthogonal indecomposable $A$-modules. Finally, we exhibit a purely representation-theoretic map to $\mathcal{A}$, providing a Caldero-Chapoton map in this setting.

\end{abstract}

\maketitle

{
\tableofcontents
}

\section*{Introduction}
Cluster algebras, introduced by Fomin and Zelevinsky in their seminal work \cite{fomin2002cluster}, are commutative algebras with a rich combinatorial structure. More precisely, a \emph{cluster algebra} is a subalgebra of a field of rational functions in $n$ variables generated by \emph{cluster variables}. Cluster variables are constructed recursively from an initial seed by a process of \emph{mutation}, and they are grouped into overlapping sets of constant cardinality $n$, called \emph{clusters}. By the \emph{Laurent phenomenon} \cite{fomin2002cluster}, every cluster variable $x$ is a Laurent polynomial in the cluster variables $u_1,\dots,u_n$ of the initial cluster. This Laurent polynomial is referred to as the \emph{cluster expansion} of $x$ in the initial cluster variables.

A cluster algebra is determined by the \emph{initial exchange matrix} $B$ and the choice of a \emph{coefficient vector} $\mathbf{y}$. A canonical choice in this setting is the \emph{principal coefficient system}, introduced in \cite{fomin2007cluster}. With this choice, $x$ turns out to be a homogeneous Laurent polynomial in the initial cluster variables. In this context, knowing the cluster expansion of $x$ is equivalent to knowing the $F$-polynomial $F_x$ and the $\mathbf{g}$-vector $\mathbf{g}_x$ of $x$, which are defined as the evaluation of $x$ at $u_1=\cdots=u_n=1$ and the multi-degree of $x$, respectively. Moreover, in \cite{fomin2007cluster}, the authors show that knowing the expansion formulas for the principal coefficients suffices to determine the expansion formulas for arbitrary coefficients.

Fomin, Shapiro and Thurston in \cite{fomin2008cluster,FT}, building on work of Fock and Goncharov \cite{FG1,FG2}, initiate the study of skew-symmetric cluster algebras $\mathcal{A}_\bullet(T)$ from triangulations $T$ of surfaces $(\textbf{S},\textbf{M})$ with boundary and marked points. In their approach, cluster variables $x_\gamma$ correspond to arcs $\gamma$ in the surface, and clusters correspond to triangulations. Then, Musiker, Schiffler and Williams in \cite{MS,MSW2011} give an expansion formula for the cluster variables in terms of perfect matchings of some labeled planar graphs, called \emph{snake graphs}, recursively constructed from the surface. Later, Felikson, Shapiro and Tumarkin in \cite{felikson2012cluster}, define skew-symmetrizable cluster algebras from triangulated orbifolds, and extend to this setting the correspondence between cluster variables and arcs, and between clusters and triangulations.

\vspace{0.5cm}

In this paper, we investigate skew-symmetrizable cluster algebras associated with triangulated surfaces equipped with an orientation-preserving $\mathbb{Z}_2$-action. Specifically, let $(\tilde{\mathbf{S}},\tilde{\mathbf{M}})$ be an unpunctured surface with an orientation-preserving diffeomorphism $\varsigma$ of order 2 that globally fixes the set of marked points $\tilde{\mathbf{M}}$. We define cluster algebras $\mathcal{A}_\bullet^\varsigma(\Tilde{T})$ with principal coefficients in certain admissible $\varsigma$-invariant triangulations $\Tilde{T}$ of $(\tilde{\mathbf{S}},\tilde{\mathbf{M}})$. This construction recovers the skew-symmetrizable cluster algebra associated by Felikson, Shapiro and Tumarkin in \cite{felikson2012cluster} to the orbifold surface $\tilde{\mathbf{S}}/\varsigma$ with one orbifold point of weight 2 and principal coefficients in the orbifold triangulation $\tilde{T}/\varsigma$. We show that admissible $\varsigma$-invariant triangulations of $(\tilde{\mathbf{S}},\tilde{\mathbf{M}})$ are in bijection with the clusters of $\mathcal{A}_\bullet^\varsigma(\Tilde{T})$, and that the cluster variables $x_{[\gamma]}$ correspond to the $\varsigma$-orbits $[\gamma]$ of the arcs of $(\Tilde{\mathbf{S}},\Tilde{\mathbf{M}})$.

Furthermore, as our first main result, we establish a formula expressing some cluster variables of $\mathcal{A}^\bullet_\varsigma(\widetilde{T})$ in terms of those of a skew-symmetric cluster algebra $\mathcal{A}^\bullet(T)$, associated with a triangulation $T=\{\tau_1,\dots,\tau_n\}$ of a surface $(\mathbf{S},\mathbf{M})$. The latter is obtained from $(\widetilde{\mathbf{S}},\widetilde{\mathbf{M}})$ by collapsing one of its two symmetric parts to a triangle. The key operation relating the two cluster algebras is the \emph{restriction} (see Definition \ref{def_restriction}). For
a $\varsigma$-orbit $[\gamma]$ of $(\Tilde{\textbf{S}},\Tilde{\textbf{M}})$, we denote by $F_{[\gamma]}$ and $\mathbf{g}_{[\gamma]}$ the $F$-polynomial and the $\mathbf{g}$-vector, respectively, of the cluster variable $x_{[\gamma]}$ of $\mathcal{A}_\bullet^\varsigma(\Tilde{T})$. On the other hand, for an arc $\gamma$ of $(\textbf{S},\textbf{M})$, we denote by $F_\gamma$ and $\mathbf{g}_{\gamma}$ the $F$-polynomial and the $\mathbf{g}$-vector, respectively, of the cluster variable $x_{\gamma}$ of $\mathcal{A}_\bullet(T)$.

   \begin{teo*}[\ref{thm:formula cv}]
    Let $[\gamma] \not\subset \Tilde{T}$ be a non-crossing $\varsigma$-orbit of $(\Tilde{\textbf{S}},\Tilde{\textbf{M}})$. Let $D=\operatorname{diag}(1,\dots,1,2)$ be the $n\times n$ diagonal matrix with diagonal entries $1,\dots,1,2$. Then,
    \begin{itemize}
        \item [(i)] If $\operatorname{Res}([\gamma])=\{\gamma_1\}$, then
        \begin{align*}
            F_{[\gamma]}&=F_{\gamma_1},\\
            \mathbf{g}_{[\gamma]}&=\begin{cases}
          \text{$D\mathbf{g}_{\gamma_1}$ \hspace{1cm}if $\gamma_1$ does not cross $\tau_n$;}\\
          \text{$D\mathbf{g}_{\gamma_1}+\mathbf{e}_n$ \hspace{0.1cm} if $\gamma_1$ crosses $\tau_n$},
      \end{cases}
        \end{align*}
        where $\tau_n$ is the unique $\varsigma$-invariant arc of $\Tilde{T}$.
        \item [(ii)] Otherwise, $\operatorname{Res}([\gamma])=\{\gamma_1,\gamma_2\}$. If $\gamma_1$ and $\gamma_2$ cross at the endpoint $\blacksquare$ in $(\textbf{S},\textbf{M})$, then
        \begin{align*}
            F_{[\gamma]}&=F_{\gamma_1}F_{\gamma_2}- \mathbf{y}^{\mathbf{d}_{\gamma_1,\gamma_2}}F_{\gamma_3},\\
            \mathbf{g}_{[\gamma]}&=D(\mathbf{g}_{\gamma_1}+\mathbf{g}_{\gamma_2}+\mathbf{e}_n),
        \end{align*}
        where $\gamma_3$ is the smoothing of the crossing of $\gamma_1$ and $\gamma_2$ at $\blacksquare$, $\mathbf{d}_{\gamma_1,\gamma_2}$ is the integer vector that keeps track of the elementary laminations of the arcs of $T$ that cross both $\gamma_1$ and $\gamma_2$, and $\mathbf{e}_n$ is the $n$-th vector of the canonical basis of $\mathbb{Z}^n$.
    \end{itemize}
\end{teo*}

Moreover, we associate with these $\varsigma$-orbits $[\gamma]$ of $(\Tilde{\textbf{S}},\Tilde{\textbf{M}})$, that we call \emph{nice}, a labeled modified snake graph $\mathcal{G}_{[\gamma]}$ constructed by gluing together the snake graphs corresponding to the arcs of $\operatorname{Res}([\gamma])$. This allows us to obtain the cluster expansion of the cluster variable $x_{[\gamma]}$ of $\mathcal{A}_\bullet^\varsigma(\Tilde{T})$ in terms of perfect matchings of $\mathcal{G}_{[\gamma]}$. Our construction generalizes to arbitrary unpunctured surfaces a previous result of the author \cite{ciliberti2}, which provides cluster expansions in terms of perfect matchings of modified snake graphs for cluster algebras of type $B$ and $C$, associated with regular polygons with an even number of vertices. This is the second main result of the paper:

\begin{teo*}[\ref{thm:sg}]
 Let $\Tilde{T}=\{\tau_i\}_{i=1}^{n-1} \sqcup \{\tau_n\} \sqcup \{\tau_i'\}_{i=1}^{n-1}$ be an admissible $\varsigma$-invariant triangulation of $(\Tilde{\textbf{S}},\Tilde{\textbf{M}})$. Assume further that $\tau_n$ and $\tau_{n-1}$ are edges of a triangle of $\Tilde{T}$ whose third edge is a boundary arc. Let $\mathcal{A}_\bullet^\varsigma(\Tilde{T})$ be the skew-symmetrizable cluster algebra with principal coefficients in $\Tilde{T}$. Let $[\gamma]$ be a nice $\varsigma$-orbit. Then 
 \begin{align*}
     F_{[\gamma]}=F_{\mathcal{G}_{[\gamma]}},
 \end{align*}
 and 
 \begin{align*}
     \mathbf{g}_{[\gamma]}=\mathbf{g}_{\mathcal{G}_{[\gamma]}},
 \end{align*} 
 where $F_{\mathcal{G}_{[\gamma]}}$ is the perfect matching polynomial of $\mathcal{G}_{[\gamma]}$, and $\mathbf{g}_{\mathcal{G}_{[\gamma]}}$ its $\mathbf{g}$-vector.  
\end{teo*}

\vspace{0.5cm}

In parallel, the representation theory of symmetric quivers has been developed by Derksen and Weyman in \cite{DW}, as well as Boos and Cerulli Irelli in \cite{BCI}. A \emph{symmetric quiver algebra} is a finite-dimensional algebra $A=kQ/I$ with an involution $\varrho$ of vertices and arrows that reverses the orientation of arrows, and preserves $I$. A \emph{symmetric module} over a symmetric algebra $A$ is an ordinary $A$-module equipped with some extra data that forces each dual pair of arrows of $Q$ to act anti-adjointly. Symmetric modules are of two types: \emph{orthogonal} and \emph{symplectic}. They form an additive category which is not abelian (see Section \ref{symmetric_quivers}).

In the last section, we introduce a categorification of skew-symmetrizable cluster algebras from surfaces endowed with an orientation-preserving involution via symmetric quivers. Namely, given a skew-symmetrizable cluster algebra $\mathcal{A}_\bullet^\varsigma(\Tilde{T})$, we associate a symmetric quiver algebra $A$ with it, in such a way that the non-initial cluster variables $x_N$ of $\mathcal{A}_\bullet^\varsigma(\Tilde{T})$ correspond to \emph{non-crossing} orthogonal indecomposable $A$-modules $N$. Moreover, Theorem \ref{thm:formula cv} enables us to define a Caldero-Chapoton-like map (see \cite{CC}) to $\mathcal{A}_\bullet^\varsigma(\Tilde{T})$.

For a non-crossing orthogonal indecomposable $A$-module $N$, $F_N$ and $\mathbf{g}_N$ denote the $F$-polynomial and the $\mathbf{g}$-vector, respectively, of $x_N$. On the other hand, $F_{\operatorname{Res}(N)}$ and $\mathbf{g}_{\operatorname{Res}(N)}$ are the $F$-polynomial and the $\mathbf{g}$-vector of the $A$-module $\operatorname{Res}(N)$, obtained from $N$ by assigning the trivial vector space to all vertices with index greater than $n$ (see Definition \ref{def_res}). The following theorem, which constitutes the third main result of this work, provides a purely representation-theoretic formula to compute $F_N$ and $\mathbf{g}_N$:
\begin{teo*}[\ref{cat_interpr}]
Let $N$ be a non-crossing orthogonal indecomposable $A$-module. Let $D=\operatorname{diag}(1,\dots,1,2)\in \mathbb{Z}^{n\times n}$.
\begin{itemize}
    \item [(i)] If $\operatorname{Res}(N)=(V_i,\phi_a)_{i=1}^n$ is indecomposable as $A$-module, then 
    \begin{equation*}
        \text{$F_N=F_{\operatorname{Res}(N)}$,}
    \end{equation*}
    and
    \begin{equation*}
 \mathbf{g}_{N}=\begin{cases}
          \text{$D \mathbf{g}_{\operatorname{Res}(N)}$ \hspace{1.8cm}if $\operatorname{dim} V_n =0$;}\\
          \text{$D \mathbf{g}_{\operatorname{Res}(N)}+\mathbf{e}_n$ \hspace{1cm}if $\operatorname{dim} V_n \neq 0$.}
      \end{cases}
      \end{equation*}
\item [(ii)]      Otherwise, $N=L\oplus \nabla L$. If $\operatorname{dim} \operatorname{Ext}^1(\nabla L, L)=1$, then there exists a non-split short exact sequence
    \begin{equation*}
        0 \to L \to G_1 \oplus G_2 \to \nabla L \to 0,
    \end{equation*}
    where $G_1$ and $G_2$ are orthogonal indecomposable $A$-modules of type I. Furthermore, denoting by $\overline{L}$ the kernel of a non-trivial map $L \to \tau \nabla L$ which does not factor through injective $A$-modules, and by $\underline{\nabla L}$ the image of a non-trivial map $\tau^{-1}L \to \nabla L$ which does not factor through projective $A$-modules, 
    \begin{equation*}
        F_N= F_{\operatorname{Res}(N)} - \mathbf{y}^{\operatorname{Res}(\textbf{dim}\underline{\nabla L})}F_{\operatorname{Res}(M)},
    \end{equation*}
    and
    \begin{equation*}
        \mathbf{g}_{N}=D(\mathbf{g}_{\operatorname{Res}(N)}+\mathbf{e}_n),
    \end{equation*}
where $M$ is the $\leq_{\mathrm{Ext}}$-minimum extension in $A$ between $\nabla L/\underline{\nabla L}$ and $\overline{L}$.
 \end{itemize}  
     
\end{teo*}

\vspace{0.5cm}

Several other works in the literature use different techniques to study skew-symmetrizable cluster algebras from surfaces. In \cite{canakciTumarkin}, \text{\c{C}anak\c{c}\i} and Tumarkin define snake and band graphs associated with arcs on a triangulated orbifold surface with orbifold points of weight $\frac{1}{2}$. In \cite{FST_finite_type}, Felikson, Shapiro and Tumarkin investigate a relation between skew-symmetric and skew-symmetrizable cluster algebras of finite mutation type via folding. In \cite{BanaianKelley}, Banaian and Kelley extend the construction of snake graphs to generalized cluster algebras arising from unpunctured orbifolds. Furthermore, other categorifications of skew-symmetrizable cluster algebras include the work of Geiss, Leclerc, and Schröer \cite{GLS1}, which uses locally free modules over certain Iwanaga-Gorenstein algebras, the species with potential approach by Geuenich and Labardini-Fragoso \cite{GeuLF1,GeuLF2}, and Demonet's construction \cite{Demonet}, which employs exact stably 2-Calabi-Yau categories endowed with a finite group action. On the other hand, in \cite{baziermatte2023markednonorientablesurfacescluster}, Bazier-Matte, Chan and Wright use symmetric modules to give a categorification of quasi-cluster algebras from non-orientable surfaces.

\vspace{0.5cm}

The paper is structured as follows. In Section \ref{s1}, we recall the definition of the skew-symmetric cluster algebra associated with a triangulation of an unpunctured surface. We then define the skew-symmetrizable cluster algebra associated with a surface equipped with an orientation-preserving involution $\varsigma$, showing that clusters are in bijection with admissible $\varsigma$-invariant triangulations, and that cluster variables correspond to non-crossing $\varsigma$-orbits of arcs. We conclude the section with the proof of Theorem \ref{thm:formula cv}. Section \ref{s2} begins with a brief overview of snake graphs arising from arcs on surfaces. Then, we present the construction of modified snake graphs, and prove Theorem \ref{thm:sg}. Finally, in Section \ref{s3}, we introduce the categorification via symmetric quivers of the skew-symmetrizable cluster algebras defined in Section \ref{s1}. In particular, after a recollection on symmetric representation theory, we associate a symmetric algebra $A$ to any of these cluster algebras such that non-crossing orthogonal indecomposable $A$-modules correspond to non-initial cluster variables. This leads to the proof of Theorem \ref{cat_interpr}.

\section{Cluster algebras from surfaces}\label{s1}

In this section, we first recall the definition of skew-symmetric cluster algebra from an unpunctured marked surface following \cite{fomin2008cluster}, and then introduce the notion of skew-symmetrizable cluster algebra from a surface with a $\bbZ_2$-action.

\vspace{0.2cm}

We work in the following setting:
\begin{itemize}
    \item $\textbf{S}$ is a connected oriented 2-dimensional Riemann surface with non-empty boundary $\partial\textbf{S}$;
    \item $\textbf{M} \subset \partial\textbf{S}$ is a finite set of marked points on the boundary of $\textbf{S}$ such that each connected component of $\partial\textbf{S}$ has at least one marked point on it.
\end{itemize}

Up to homeomorphism, the \emph{surface} $(\textbf{S},\textbf{M})$ is determined by:
\begin{itemize}
    \item the genus of $\textbf{S}$;
    \item the number of boundary components;
    \item the number of marked points on each boundary component.
\end{itemize}

\subsection{Skew-symmetric cluster algebras from surfaces}

\begin{definition}[Arc]\label{def:arcs}
    An \emph{arc} $\gamma$ in $(\textbf{S},\textbf{M})$ is a curve in $\textbf{S}$ such that
    \begin{itemize}
        \item [(a)] the endpoints of $\gamma$ are in $\textbf{M}$;
        \item [(b)] $\gamma$ does not intersect itself, except that its endpoints may coincide;
        \item [(c)] except for the endpoints, $\gamma$ is disjoint from $\textbf{M}$ and $\partial \textbf{S}$;
        \item [(d)] $\gamma$ is not contractible into $\textbf{M}$ or into $\partial\textbf{S}$.
    \end{itemize}
    Curves which satisfy (a), (b), (c), and are contractible into $\partial\textbf{S}$, are called \emph{boundary arcs}. Each arc or boundary arc is considered up to isotopy.
\end{definition}

\begin{definition}[Compatible arcs]\label{def:compatible}
    Two arcs are called \emph{compatible} if they do not intersect in the interior of $\textbf{S}$; more precisely, there are curves in their respective isotopy classes which do not intersect in the interior of $\textbf{S}$.
\end{definition}

\begin{definition}[Triangulation]\label{def:ideal triang}
    A maximal collection of distinct pairwise compatible arcs is called a \emph{triangulation}. The arcs of a triangulation cut the surface $\textbf{S}$ into \emph{triangles}.
\end{definition}

\begin{remark}
    An elementary topological argument shows that the number $n$ of arcs in a triangulation is an invariant of $(\textbf{S},\textbf{M})$, known as the \emph{rank} of $(\textbf{S},\textbf{M})$.
\end{remark}

\begin{definition}[Flip]\label{def:flip}
    A \emph{flip} is a transformation of a triangulation $T$ that removes an arc $\gamma$ and replaces it with the unique arc $\gamma' \neq \gamma$ that, together with the remaining arcs, forms a new triangulation $T'$.
\end{definition}

All triangulations of $(\textbf{S},\textbf{M})$ are connected by a series of flips.

\begin{definition}[Signed adjacency matrix]\label{def:adj_matrix}
    To each triangulation $T=\{ \tau_1,\dots, \tau_n \}$ is associated the \emph{signed adjacency matrix} $B(T)$ defined in the following way:
    \begin{itemize}
        \item for each triangle $\Delta$ in $T$, consider the $n \times n$ integer matrix $B^\Delta=(b_{ij}^\Delta)$, where
\begin{equation*}
b_{ij}^\Delta=
\begin{cases}
1 & \!\!\text{if $\tau_i$ and $\tau_j$ are sides of $\Delta$,}\\
& \hspace{-0.5cm}\text{\quad with $\tau_i$ following $\tau_j$ in  counterclockwise order;}\\
-1 & \!\!\text{if $\tau_i$ and $\tau_j$ are sides of $\Delta$,}\\
& \hspace{-0.5cm}\text{\quad with $\tau_j$ following $\tau_i$ in  counterclockwise order;}\\
0 & \!\!\text{otherwise;}
\end{cases}
\end{equation*}
\item the matrix $B=B(T)=(b_{ij})$ is then defined by
\begin{center}
    $B=\displaystyle\sum_\Delta B^\Delta$,
\end{center}
where the sum is over all triangles $\Delta$ in $T$.
    \end{itemize}
\end{definition} 

\begin{definition}[Cluster algebra with principal coefficients in $T$]\label{def:ca_triang}
    Let $T$ be a triangulation of $(\textbf{S},\textbf{M})$. The \emph{cluster algebra $\mathcal{A}_\bullet(T)$ associated to the surface $(\textbf{S},\textbf{M})$ with principal coefficients in $T$} is defined as the cluster algebra with principal coefficients in the initial seed whose exchange matrix is $B(T)$.
\end{definition}

Fomin, Shapiro and Thurston proved the following correspondence:

\begin{theorem}[\cite{fomin2008cluster}]\label{thm:fst}
There are bijections:
   \[
\begin{array}{ccc}
\big\{\text{cluster variables of } \mathcal{A}_\bullet(T)\big\} & \longleftrightarrow & \big\{\text{arcs of } (\textbf{S},\textbf{M})\big\} \\
x_\gamma & & \gamma \\[15pt]
\big\{\text{clusters of } \mathcal{A}_\bullet(T)\big\} & \longleftrightarrow & \big\{\text{triangulations of } (\textbf{S},\textbf{M})\big\} \\
\mathbf{x}_{\mathcal{T}} = \{x_{\tau_1}, \dots, x_{\tau_n}\} & & \mathcal{T} = \{\tau_1, \dots, \tau_n\}
\end{array}
\] 
\end{theorem}

\begin{remark}
   The initial cluster of $\mathcal{A}_\bullet(T)$ is $\mathbf{x}_{T}$. We denote the initial cluster variables by $u_1,\dots,u_n$, where $u_i$ is cluster variable corresponding to the arc $\tau_i$. Moreover, if $\gamma$ is a boundary arc, then $x_{\gamma}=1$. Finally, we denote by $F_{\gamma}$ and $\mathbf{g}_{\gamma}$ the $F$-polynomial and the $\mathbf{g}$-vector of $x_{\gamma}$, respectively. 
\end{remark}

Furthermore, the exchange relations in $\mathcal{A}_\bullet(T)$ correspond to flips of the arcs, and the coefficients are given by elementary laminations of the arcs of $T$.

\begin{definition}[Elementary lamination of $\gamma$] \label{def:elem-laminate}
Let $\gamma$ be an arc of $(\textbf{S},\textbf{M})$. The \emph{elementary lamination} associated to $\gamma$ is the arc $L_{\gamma}$ which runs along $\gamma$ within a small neighborhood of it. In particular, if $\gamma$ begins at a marked point $a$ on a connected component $C_1$ of $\partial\textbf{S}$ and ends at a marked point $b$ on a connected component $C_2$ of $\partial\textbf{S}$, then $L_\gamma$ begins at a point $a' \in C_1$ located near $a$ in
the clockwise direction and ends at a point $b' \in C_2$ near $b$ in
the clockwise direction. If $T=\{ \tau_1,\dots,\tau_n \}$ is a triangulation of $(\textbf{S},\textbf{M})$, then $L_{\tau_i}$ is denoted by $L_i$.
\end{definition}

\begin{notation} We use the following notation:
\begin{itemize}
    \item [(i)] Let $\mathbf{d}=(d_i) \in \mathbb{Z}_{\geq 0}^n$ be an integer vector. We denote by $\mathbf{y}^{\mathbf{d}}$ the monomial $y_1^{d_1}\cdots y_n^{d_n}$.
    \item [(ii)] Let $T=\{\tau_1,\dots,\tau_n\}$ be a triangulation of $(\textbf{S},\textbf{M})$. Given two arcs $\gamma_1$ and $\gamma_2$ of $(\textbf{S},\textbf{M})$, we denote by $\mathbf{d}_{\gamma_1,\gamma_2}$ the integer vector $\mathbf{d}_{\gamma_1,\gamma_2}=(d_i)\in \mathbb{Z}_{\geq 0}^n$ whose $i$-th coordinate $d_i$ is given by the number of times that $L_{\tau_i}$ crosses both $\gamma_1$ and $\gamma_2$.
\end{itemize} 
\end{notation}

\begin{definition}[Smoothing of a crossing of two arcs at an interior point]
    Let $\gamma_1$ and $\gamma_2$ be two arcs that cross at an interior point $x$. The \emph{smoothing of the crossing of $\gamma_1$ and $\gamma_2$ at the point $x$} is given by the pairs of arcs $\{\alpha_1, \alpha_2\}$ and $\{\beta_1, \beta_2\}$ such that

\begin{itemize}
    \item[-] $\{\alpha_1, \alpha_2\}$ is the same as $\{\gamma_1, \gamma_2\}$ except locally where the crossing $\times$ is replaced by the pair of segments $\asymp$,
    \item[-] $\{\beta_1, \beta_2\}$ is the same as $\{\gamma_1, \gamma_2\}$ except locally where the crossing $\times$ is replaced by the pair of segments $\supset \subset$.
\end{itemize}
The roles of the two pairs $\{\alpha_1,\alpha_2\}$ and $\{\beta_1,\beta_2\}$ are interchangeable.
\end{definition}

\begin{definition}[Smoothing of a crossing of two arcs at an endpoint]
    Let $\gamma_1$ and $\gamma_2$ be two arcs that cross at an endpoint $x$. The \emph{smoothing of the crossing of $\gamma_1$ and $\gamma_2$ at the point $x$} is the arc $\alpha$ such that $\alpha$ is the same as the concatenation of $\gamma_1$ and $\gamma_2$ except locally where the crossing $\times$ is replaced by the segment $\supset$.
\end{definition}

\begin{proposition}[\cite{MW13}] \label{up:skein1}
Let $\mathcal{A}_\bullet(T)$ be the cluster algebra associated to $(\textbf{S},\textbf{M})$ with principal coefficients in the triangulation $T=\{ \tau_1, \dots, \tau_n \}$.
Let $\gamma_1$ and $\gamma_2$ be two arcs
that cross at an interior point $x$, and let $\{\alpha_1, \alpha_2\}$ and $\{\beta_1, \beta_2\}$ be the smoothing of the crossing of $\gamma_1$ and $\gamma_2$ at $x$. Then
\begin{equation*} \label{u:skein-eq1}
x_{\gamma_1} x_{\gamma_2} = \mathbf{y}^{\mathbf{d}_{\beta_1,\beta_2}} x_{\alpha_1} ~x_{\alpha_2} + \mathbf{y}^{\mathbf{d}_{\alpha_1,\alpha_2}}x_{\beta_1} ~x_{\beta_2}.
\end{equation*}
\end{proposition}

\subsection{Skew-symmetrizable cluster algebras from surfaces with a $\bbZ_2$-action}\label{section:ssca_surfaces}

In this section, we focus on triangulated surfaces endowed with an orientation-preserving $\bbZ_2$-action, and define the skew-symmetrizable cluster algebras associated with them, although the first two definitions do not strictly require the action to be orientation-preserving.

\begin{figure}[h]
    \centering
    \includegraphics[width=0.3\linewidth]{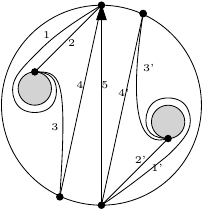}
    \caption{Example of admissible $\varsigma$-invariant triangulation.}
    \label{fig:ex_adm_triang}
\end{figure}

\begin{definition}[Admissible $\varsigma$-invariant triangulation]
Let $(\Tilde{\textbf{S}},\Tilde{\textbf{M}})$ be a surface endowed with a diffeomorphism $\varsigma$ of order 2. Let $\Tilde{T}$ be a $\varsigma$-invariant triangulation of $(\Tilde{\textbf{S}},\Tilde{\textbf{M}})$. We say that $\Tilde{T}$ is \emph{admissible} if it satisfies the following properties:
    \begin{itemize}
    \item [(i)] $\Tilde{T}$ has exactly one $\varsigma$-invariant arc;
    \item [(ii)] there is a relabeling of the arcs of $\Tilde{T}$ that induces a partition $\Tilde{T}=\{\tau_i\}_{i=1}^{n-1} \sqcup \{\tau_n\} \sqcup \{\tau_i'\}_{i=1}^{n-1}$ such that 
    \begin{itemize}
        \item $\varsigma(\tau_i)=\tau_i'$ for any $i=1,\dots,n-1$;
        \item $\varsigma(\tau_n)=\tau_n$;
        \item for any $i,j=1,\dots,n-1$, no triangle $\Delta$ of $\Tilde{T}$ has both $\tau_i$ and $\tau_j'$ as edges.
    \end{itemize}
\end{itemize}
Furthermore, we assume that $\tau_n$ is oriented. See Figure \ref{fig:ex_adm_triang} for an example.
\end{definition}

\begin{example}
   The triangulation in Figure \ref{fig:non_adm} is not admissible.
    \begin{figure}[h]
    \centering
    \includegraphics[width=0.3\linewidth]{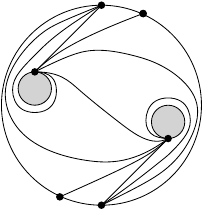}
    \caption{Example of non-admissible $\varsigma$-invariant triangulation.}
    \label{fig:non_adm}
\end{figure}
\end{example}

\begin{definition}[Restriction]\label{def_restriction}
Let $(\Tilde{\textbf{S}},\Tilde{\textbf{M}})$ be a surface endowed with a diffeomorphism $\varsigma$ of order 2. Let $\mathcal{D}$ be a set of arcs of $(\Tilde{\textbf{S}},\Tilde{\textbf{M}})$, and let $\Tilde{T}$ be an admissible $\varsigma$-invariant triangulation. We define the \emph{restriction of $\mathcal{D}$}, and we denote it by $\operatorname{Res}(\mathcal{D})$, as the set of arcs of the collapsed surface $(\textbf{S},\textbf{M})$ resulting from the collapse of the region to the right of $\tau_n$ to a triangle. 
\end{definition}
The marked point of $(\textbf{S},\textbf{M})$ resulting from the collapse of the region to the right of $\tau_n$ is denoted by $\blacksquare$.

\begin{figure}[h]
    \centering
    \includegraphics[width=0.85\linewidth]{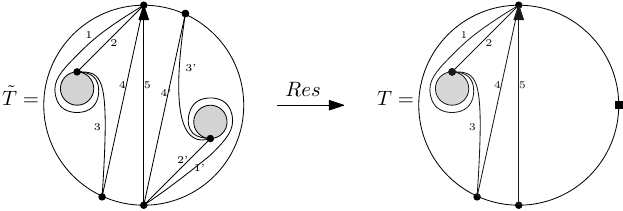}
    \caption{On the left, an admissible $\varsigma$-invariant triangulation of a surface of genus 0 with three boundary components; on the right, its restriction.}
    \label{fig:restriction}
\end{figure}

\begin{remark}
    The collapsed surface $(\textbf{S},\textbf{M})$ is a surface of rank $n$. In particular, $T=\operatorname{Res}(\Tilde{T})=\{\tau_1, \dots, \tau_n\}$ is a triangulation of $(\textbf{S},\textbf{M})$. See Figure \ref{fig:restriction} for an example.
\end{remark}

In the following, $(\Tilde{\textbf{S}},\Tilde{\textbf{M}})$ is a surface together with an orientation-preserving diffeomorphism $\varsigma$ of order 2, fixing globally $\Tilde{\textbf{M}}$. 

\begin{definition}[Skew-symmetrizable cluster algebra with principal coefficients in $\Tilde{T}$]\label{def:ss_ca_tilde_T}
 Let $\Tilde{T}$ be an admissible $\varsigma$-invariant triangulation of $(\Tilde{\textbf{S}},\Tilde{\textbf{M}})$. We define the \emph{cluster algebra $\mathcal{A}_\bullet^\varsigma(\Tilde{T})$ associated to the surface $(\Tilde{\textbf{S}},\Tilde{\textbf{M}})$ with principal coefficients in $\Tilde{T}$} as the cluster algebra with principal coefficients in the initial seed whose exchange matrix is $B^\varsigma(\Tilde{T}):=DB(\operatorname{Res}(\Tilde{T}))$, where $D=\operatorname{diag}(1,\dots,1,2)$ is the $n\times n$ diagonal matrix with diagonal entries $1,\dots,1,2$.   
\end{definition}

\begin{remark}
    The cluster algebra $\mathcal{A}_\bullet^\varsigma(\Tilde{T})$ just defined is the skew-symmetrizable cluster algebra associated in \cite{felikson2012cluster} with the orbifold surface $\Tilde{\textbf{S}}/\varsigma$ with one orbifold point of weight 2 corresponding to the $\varsigma$-invariant arc $\tau_n$.
\end{remark}

\begin{remark}
    If $\tilde{\mathbf{S}}$ is a regular polygon with $2n+2$ vertices, any  triangulation invariant under 180°-rotation is admissible, and $\mathcal{A}_\bullet^\varsigma(\Tilde{T})$ is a cluster algebra of type $B_n$.
\end{remark}

Let $\gamma$ be an arc of $(\Tilde{\textbf{S}},\Tilde{\textbf{M}})$. We denote by $[\gamma]$ the $\varsigma$-orbit of $\gamma$. A $\varsigma$-orbit can be either a $\varsigma$-invariant arc or a $\varsigma$-invariant pair of non-$\varsigma$-invariant arcs. We say that a $\varsigma$-orbit is \emph{non-crossing} if its arcs do not cross in the interior. Note that if $\Tilde{\textbf{S}}$ is a regular polygon, then every $\varsigma$-orbit is non-crossing.

\begin{proposition}\label{prop:cluster var in A^sigma}
Let $\Tilde{T}$ be an admissible $\varsigma$-invariant triangulation of a surface $(\Tilde{\textbf{S}},\Tilde{\textbf{M}})$ endowed with an orientation-preserving diffeomorphism $\varsigma$ of order 2, fixing globally $\Tilde{\textbf{M}}$. Then, there are bijections:
   \[
\begin{array}{ccc}
\big\{\text{cluster variables of } \mathcal{A}_\bullet^\varsigma(\Tilde{T})\big\} & \longleftrightarrow & \big\{\text{non-crossing $\varsigma$-orbits of arcs of } (\Tilde{\textbf{S}},\Tilde{\textbf{M}})\big\} \\
x_{[\gamma]} & & [\gamma] \\[15pt]
\big\{\text{clusters of } \mathcal{A}_\bullet^\varsigma(\Tilde{T})\big\} & \longleftrightarrow & \big\{\text{$\varsigma$-invariant admissible triangulations of } (\Tilde{\textbf{S}},\Tilde{\textbf{M}})\big\} \\
\mathbf{x}_{\mathcal{T}} = \{x_{[\tau_1]}, \dots, x_{[\tau_n]}\} & & \mathcal{T} = \{\tau_i\}_{i=1}^{n-1} \sqcup \{\tau_n\} \sqcup \{\tau_i'\}_{i=1}^{n-1}
\end{array}
\] 
\end{proposition}

\begin{proof}
     Admissible $\varsigma$-invariant triangulations of $(\Tilde{\textbf{S}},\Tilde{\textbf{M}})$ correspond to triangulations of the orbifold surface $\Tilde{\textbf{S}}/\varsigma$ via folding (see \cite[Section 12]{felikson2012cluster} for more details). Moreover, non-crossing $\varsigma$-orbits of arcs of $(\Tilde{\textbf{S}},\Tilde{\textbf{M}})$ correspond to arcs of $\Tilde{\textbf{S}}/\varsigma$. Consequently, the statement is a reformulation of \cite[Theorem 5.9]{felikson2012cluster}.
\end{proof}

\begin{remark}
    The initial cluster of $\mathcal{A}_\bullet^\varsigma(\Tilde{T})$ is $\mathbf{x}_{\Tilde{T}}$. We denote the initial cluster variables by $u_1,\dots,u_n$, where $u_i$ is cluster variable corresponding to the $\varsigma$-orbit $[\tau_i]$. Moreover, if $\gamma$ is a boundary arc, then $x_{[\gamma]}=1$. Finally, we denote by $F_{[\gamma]}$ and $\mathbf{g}_{[\gamma]}$ the $F$-polynomial and the $\mathbf{g}$-vector of $x_{[\gamma]}$, respectively.
\end{remark}

\begin{remark}
    If $\mathbf{S}$ is a regular polygon, we recover \cite[Section 12.3]{fomin2003cluster}.
\end{remark}

\begin{notation}
Let $\Tilde{T}=\{\tau_i\}_{i=1}^{n-1} \sqcup \{\tau_n\} \sqcup \{\tau_i'\}_{i=1}^{n-1}$ be a $\varsigma$-invariant triangulation of $(\Tilde{\textbf{S}},\Tilde{\textbf{M}})$. Given a $\varsigma$-orbit $[\gamma]=\{\gamma,\gamma'\}$, we denote by $\mathbf{d}_{[\gamma]}$ the integer vector $\mathbf{d}_{[\gamma]}=(d_i)\in \mathbb{Z}_{\geq 0}^n$ whose $i$-th coordinate is given by
    \begin{align*}
            d_i=(\mathbf{d}_{\gamma,\gamma'})_i.
        \end{align*}
Similarly, given two $\varsigma$-orbits $[\gamma]=\{\gamma,\gamma'\}$ and $[\delta]=\{\delta,\delta'\}$, with $\gamma \neq \gamma'$, we denote by $\mathbf{d}_{[\gamma],[\delta]}$ the integer vector $\mathbf{d}_{[\gamma],[\delta]}=(d_i)\in \mathbb{Z}_{\geq 0}^n$ whose $i$-th coordinate is given by
    \begin{align*}
            d_i=(\mathbf{d}_{\gamma,\delta})_i+(\mathbf{d}_{\gamma',\delta'})_i.
        \end{align*}

\end{notation}

\begin{figure}[h]
    \centering
    \includegraphics[width=1\linewidth]{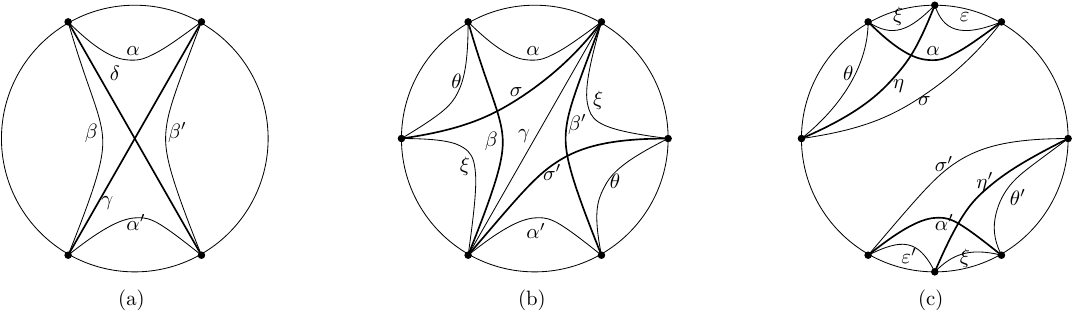}
    \caption{Skein relations in $\mathcal{A}_\bullet^\varsigma(\Tilde{T})$.}
    \label{fig:skein_a_tilde}
\end{figure}

\begin{proposition}\label{prop: relations in A^sigma}
   The following relations hold in $\mathcal{A}_\bullet^\varsigma(\Tilde{T})$:
    \begin{itemize}
        \item [(a)] Let $\gamma$, $\delta$ be two $\varsigma$-invariant arcs of $(\Tilde{\textbf{S}},\Tilde{\textbf{M}})$ as in Figure \ref{fig:skein_a_tilde}(a). Then, in the notation of Figure \ref{fig:skein_a_tilde}(a),
        \begin{equation*}
        x_{[\gamma]}x_{[\delta]}=\mathbf{y}^{\mathbf{d}_{[\beta]}}x_{[\alpha]}+\mathbf{y}^{\mathbf{d}_{[\alpha]}}x_{[\beta]}.
        \end{equation*}
        \item [(b)] Let $\beta, \sigma$ be two non-$\varsigma$-invariant arcs as in Figure \ref{fig:skein_a_tilde}(b), where $\gamma$ is a $\varsigma$-invariant arc. Then, in the notation of Figure \ref{fig:skein_a_tilde}(b),
        \begin{equation*}
        x_{[\beta]}x_{[\sigma]}=\mathbf{y}^{\mathbf{d}_{[\alpha],[\xi]}}x_{[\gamma]}^2x_{[\theta]}+\mathbf{y}^{\mathbf{d}_{[\gamma],[\theta]}}x_{[\alpha]}x_{[\xi]}.
        \end{equation*}
        \item [(c)] Let $\alpha, \eta$ be two non-$\varsigma$-invariant arcs as in Figure \ref{fig:skein_a_tilde}(c). Then, in the notation of Figure \ref{fig:skein_a_tilde}(c),
        \begin{equation*}
        x_{[\alpha]}x_{[\eta]}=\mathbf{y}^{\mathbf{d}_{[\varepsilon],[\theta]}}x_{[\sigma]}x_{[\zeta]}+\mathbf{y}^{\mathbf{d}_{[\sigma],[\zeta]}}x_{[\varepsilon]}x_{[\theta]}.
        \end{equation*}
    \end{itemize}

\end{proposition}

\begin{proof}
    The statement follows immediately from the description of the skein relations on the orbifold surface $\Tilde{\mathbf{S}}/\varsigma$, with one orbifold point of weight 2, with respect to the triangulation $\Tilde{T}/\varsigma$. See \cite[Lemma 5.8]{felikson2012cluster} and \cite[Lemma 5.26]{FelTum} for further details.
\end{proof}

\begin{theorem}\label{thm:formula cv}
    Let $[\gamma] \not\subset \Tilde{T}$ be a non-crossing $\varsigma$-orbit of $(\Tilde{\textbf{S}},\Tilde{\textbf{M}})$. Let $D=\operatorname{diag}(1,\dots,1,2)$ be the $n\times n$ diagonal matrix with diagonal entries $1,\dots,1,2$. Then,
    \begin{itemize}
        \item [(i)] If $\operatorname{Res}([\gamma])=\{\gamma_1\}$, then
        \begin{align*}
            F_{[\gamma]}&=F_{\gamma_1},\\
            \mathbf{g}_{[\gamma]}&=\begin{cases}
          \text{$D\mathbf{g}_{\gamma_1}$ \hspace{1cm}if $\gamma_1$ does not cross $\tau_n$;}\\
          \text{$D\mathbf{g}_{\gamma_1}+\mathbf{e}_n$ \hspace{0.1cm} if $\gamma_1$ crosses $\tau_n$}.
      \end{cases}
        \end{align*}
        \item [(ii)] Otherwise, $\operatorname{Res}([\gamma])=\{\gamma_1,\gamma_2\}$. If $\gamma_1$ and $\gamma_2$ cross at the endpoint $\blacksquare$ in $(\textbf{S},\textbf{M})$, then
        \begin{align*}
            F_{[\gamma]}&=F_{\gamma_1}F_{\gamma_2}- \mathbf{y}^{\mathbf{d}_{\gamma_1,\gamma_2}}F_{\gamma_3},\\
            \mathbf{g}_{[\gamma]}&=D(\mathbf{g}_{\gamma_1}+\mathbf{g}_{\gamma_2}+\mathbf{e}_n),
        \end{align*}
        where $\gamma_3$ is the smoothing of the crossing of $\gamma_1$ and $\gamma_2$ at $\blacksquare$, and $\mathbf{e}_n$ is the $n$-th vector of the canonical basis of $\mathbb{Z}^n$.
    \end{itemize}
\end{theorem}

    \begin{figure}[h] 
                               \centering
             \includegraphics[width=0.7\linewidth]{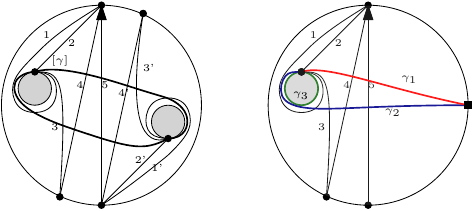}
                               \caption{A $\varsigma$-orbit $[\gamma]$ (left) and its restriction (right).}
                               \label{fig:sigma_orbit}
                               \end{figure}

\begin{example}
Consider the $\varsigma$-orbit in Figure \ref{fig:sigma_orbit},
\begin{align*}
    F_{[\gamma]}&=F_{\gamma_1}F_{\gamma_2}-y_1y_3y_4y_5F_{\gamma_3}\\
    &=(y_4y_5 + y_4 + 1)(y_1y_3y_4y_5 + y_1y_3y_4 + y_1y_4y_5 + y_1y_4 + y_4y_5 + y_1 + y_4 + 1)-y_1y_3y_4y_5\\
    &=y_1y_3y_4^2y_5^2 + 2y_1y_3y_4^2y_5 + y_1y_4^2y_5^2 + y_1y_3y_4^2 + 2y_1y_4^2y_5 + y_4^2y_5^2 + y_1y_3y_4 + y_1y_4^2 + 2y_1y_4y_5 \\&+ 2y_4^2y_5 + 2y_1y_4 + y_4^2 + 2y_4y_5 + y_1 + 2y_4 + 1;\\
    \mathbf{g}_{[\gamma]}&=D(\mathbf{g}_{\gamma_1}+\mathbf{g}_{\gamma_2}+\mathbf{e}_5)=D(\displaystyle\begin{pmatrix}
    0\\0 \\ 1 \\ -1\\0
\end{pmatrix}+\displaystyle\begin{pmatrix}
    -1 \\ 1 \\ 1\\-1\\0
\end{pmatrix}+\displaystyle\begin{pmatrix}
    0 \\0\\0\\ 0 \\ 1
\end{pmatrix})=D(\displaystyle\begin{pmatrix}
    -1\\1\\2\\-2\\1
\end{pmatrix})=\displaystyle\begin{pmatrix}
    -1\\1\\2\\-2\\2
\end{pmatrix}.
\end{align*}

\end{example}

\begin{remark}
    If $(\Tilde{\textbf{S}},\Tilde{\textbf{M}})$ is a regular polygon, we recover \cite[Theorem 3.7]{ciliberti1}. 
\end{remark}

\begin{proof}[Proof of Theorem \ref{thm:formula cv}]
   We prove the theorem by induction on the number $k\geq 1$ of intersections between each arc of $[\gamma]$ and $\Tilde{T}$. Assume $k=1$. Since each arc of $[\gamma]$ crosses only one arc of $\Tilde{T}$, either $[\gamma]$ is a pair of arcs which do not cross $\tau_n$, or $[\gamma]$ is the $\varsigma$-invariant arc that crosses only $\tau_n$ once. Therefore, $\operatorname{Res}([\gamma])=\{\gamma_j\}$, where $\gamma_j$ is the arc of $(\textbf{S},\textbf{M})$ that crosses only $\tau_j$ once. Let $T=\operatorname{Res}(\Tilde{T})$, and let $B^\varsigma(\Tilde{T})=DB(T)=(\Tilde{b}_{ij})$ and $B(T)=(b_{ij})$. 
    We have 
    \begin{equation*}
        x_{[\gamma]}u_j=y_j\displaystyle\prod_{\Tilde{b}_{ij}>0}u_i^{\Tilde{b}_{ij}}+\displaystyle\prod_{\Tilde{b}_{ij}<0}u_i^{-\Tilde{b}_{ij}}\hspace{1cm}\text{in }\mathcal{A}_\bullet^\varsigma(\Tilde{T}),
    \end{equation*}
    and
     \begin{equation*}
        x_{\gamma_j}u_j=y_j\displaystyle\prod_{b_{ij}>0}u_i^{b_{ij}}+\displaystyle\prod_{b_{ij}<0}u_i^{-b_{ij}}\hspace{1cm}\text{in }\mathcal{A}_\bullet(T).
    \end{equation*}
    So
    \begin{equation*}
        F_{[\gamma]}=y_j+1=F_{\gamma_j}.
    \end{equation*}
    If $t\neq n$,
\begin{equation*}
(\mathbf g_{[\gamma]})_t
=
\left(\deg\frac{\prod_{\widetilde b_{ij}<0}u_i^{-\widetilde b_{ij}}}{u_j}\right)_t
=
\left(\deg\frac{\prod_{b_{ij}<0}u_i^{-b_{ij}}}{u_j}\right)_t
=
(\mathbf g_{\gamma_j})_t.
\end{equation*}
If $t=n$ and $j \neq n$,
\begin{equation*}
(\mathbf g_{[\gamma]})_n
=
\left(
\deg\frac{\prod_{\widetilde b_{ij}<0}u_i^{-\widetilde b_{ij}}}{u_j}
\right)_n
=
2\left(
\deg\frac{\prod_{b_{ij}<0}u_i^{-b_{ij}}}{u_j}
\right)_n
=
2(\mathbf g_{\gamma_j})_n.
\end{equation*}
Finally, if $t=n$ and $j=n$,
\begin{equation*}
(\mathbf g_{[\gamma]})_n
=
\left(\deg\frac{1}{u_n}\right)_n
=
-1
=
(\mathbf g_{\gamma_n})_n
=
2(\mathbf g_{\gamma_n})_n+1.
\end{equation*}
Assume now $k>1$. There are three cases to consider.
    \begin{itemize}
\item [1)] Let $[\gamma]=\{\gamma,\gamma'\}$, $\gamma \neq \gamma'$, be a $\varsigma$-invariant pair of non-$\varsigma$-invariant arcs such that $\operatorname{Res}([\gamma])=\{\gamma\}$. Let $\tau_{i_1}$ be the first arc of $\Tilde{T}$ crossed by $\gamma$.
\begin{figure}[h]
    \centering
    \includegraphics[width=0.7\linewidth]{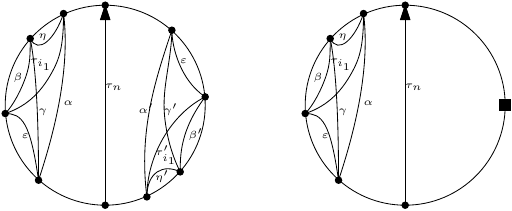}
    \caption{On the left, the skein relation between $[\gamma]$ and $[\tau_{i_1}]$; on the right, the skein relation between $\operatorname{Res}([\gamma])$ and $\operatorname{Res}([\tau_{i_1}])$ in the collapsed surface.}
    \label{fig:case1}
\end{figure}
By Proposition \ref{prop: relations in A^sigma} (c), in the notation of Figure \ref{fig:case1}, we have
\begin{align}\label{thm1_eqn1}
u_{i_1}x_{[\gamma]}&=\mathbf{y}^{\mathbf{d}_{[\varepsilon],[\eta]}}x_{[\alpha]}x_{[\beta]}+\mathbf{y}^{\mathbf{d}_{[\alpha],[\beta]}}x_{[\varepsilon]}x_{[\eta]}\notag\\
&=\mathbf{y}^{\mathbf{d}_{\varepsilon,\eta}}x_{[\alpha]}x_{[\beta]}+\mathbf{y}^{\mathbf{d}_{\alpha,\beta}}x_{[\varepsilon]}x_{[\eta]}.
\end{align}

Since $\tau_{i_1}$ is the first arc of $\Tilde{T}$ that is crossed by $\gamma$, $\beta$ and $\eta$ must be either boundary arcs or arcs of $\Tilde{T}$. It follows from \eqref{thm1_eqn1} that
\begin{equation*}
    F_{[\gamma]}=\mathbf{y}^{\mathbf{d}_{\varepsilon,\eta}}F_{[\alpha]}+\mathbf{y}^{\mathbf{d}_{\alpha,\beta}}F_{[\varepsilon]}.
\end{equation*}
Since each arc of $\Tilde{T}$ that crosses $\alpha$ (resp. $\varepsilon$) also crosses $\gamma$, the number of intersections between $\alpha$ (resp. $\varepsilon$) and $\Tilde{T}$ is strictly lower than the number of intersections between $\gamma$ and $\Tilde{T}$. By inductive hypothesis and Proposition \ref{up:skein1},
\begin{equation*}
    F_{[\gamma]}=\mathbf{y}^{\mathbf{d}_{\varepsilon,\eta}}F_{\alpha}+\mathbf{y}^{\mathbf{d}_{\alpha,\beta}}F_{\varepsilon}=F_\gamma.
\end{equation*}

   It also follows from \eqref{thm1_eqn1} that
\begin{equation*}
    \mathbf{e}_{i_1}+\mathbf{g}_{[\gamma]}=\begin{cases}
        \mathbf{g}_{[\alpha]}+\mathbf{g}_{[\beta]}\hspace{1cm} \text{if $\mathbf{y}^{\mathbf{d}_{\varepsilon,\eta}}=1$}\\
        \mathbf{g}_{[\varepsilon]}+\mathbf{g}_{[\eta]}\hspace{1cm} \text{otherwise}.
    \end{cases}
\end{equation*}
By inductive hypothesis and Proposition \ref{up:skein1},

\begin{equation*}
    \mathbf{g}_{[\gamma]}=\begin{cases}
       D(-\mathbf{e}_{i_1}+ \mathbf{g}_{\alpha}+\mathbf{g}_{\beta}) \hspace{0.5cm}\text{if $\mathbf{y}^{\mathbf{d}_{\varepsilon,\eta}}=1$}\\
    D(-\mathbf{e}_{i_1} + \mathbf{g}_{\varepsilon}+\mathbf{g}_{\eta})\hspace{0.6cm} \text{otherwise}.
        \end{cases}\hspace{-0.3cm}=D\mathbf{g}_{\gamma}.
\end{equation*} 
\item [2)] If $[\gamma]=\{\gamma\}$ is a $\varsigma$-invariant arc, the result follows by applying Proposition \ref{prop: relations in A^sigma} (a) to the crossing of $\gamma$ and $\tau_n$. The proof is omitted as it is similar to that of the first case.

\item [3)] Let $[\gamma]=\{\gamma,\gamma'\}$ be a $\varsigma$-invariant pair of non-$\varsigma$-invariant arcs such that $\operatorname{Res}([\gamma])=\{\gamma_1,\gamma_2\}$. Let $\tau_{i_1}$ be the first arc of $\Tilde{T}$ crossed by $\gamma$.
\begin{figure}[h]
    \centering
    \includegraphics[width=0.7\linewidth]{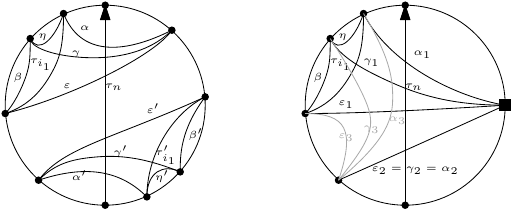}
    \caption{On the left, the skein relation between $[\gamma]$ and $[\tau_{i_1}]$; on the right, the skein relation between $\operatorname{Res}([\gamma])$ and $\operatorname{Res}([\tau_{i_1}])$ in the collapsed surface.}
    \label{fig:case3}
\end{figure}
By Proposition \ref{prop: relations in A^sigma} (c), in the notation of Figure \ref{fig:case3}, we have

\begin{align}\label{thm1_eqn1_case2}
u_{i_1}x_{[\gamma]}&=\mathbf{y}^{\mathbf{d}_{[\varepsilon],[\eta]}}x_{[\alpha]}x_{[\beta]}+\mathbf{y}^{\mathbf{d}_{[\alpha],[\beta]}}x_{[\varepsilon]}x_{[\eta]}\notag\\
&=\mathbf{y}^{\mathbf{d}_{\varepsilon_1,\eta}}x_{[\alpha]}x_{[\beta]}+\mathbf{y}^{\mathbf{d}_{\alpha_1,\beta}}x_{[\varepsilon]}x_{[\eta]}.
\end{align}
Since $\tau_{i_1}$ is the first arc of $\Tilde{T}$ that is crossed by $\gamma$, $\beta$ and $\eta$ must be either boundary arcs or arcs of $\Tilde{T}$. It follows from \eqref{thm1_eqn1_case2} that
\begin{equation*}
    F_{[\gamma]}=\mathbf{y}^{\mathbf{d}_{\varepsilon_1,\eta}}F_{[\alpha]}+\mathbf{y}^{\mathbf{d}_{\alpha_1,\beta}}F_{[\varepsilon]}.
\end{equation*}
Since each arc of $\Tilde{T}$ that crosses $\alpha$ (resp. $\varepsilon$) also crosses $\gamma$, the number of intersections between $\alpha$ (resp. $\varepsilon$) and $\Tilde{T}$ is strictly lower than the number of intersections between $\gamma$ and $\Tilde{T}$. By inductive hypothesis and Proposition \ref{up:skein1},
\begin{align*}
    F_{[\gamma]}&=\mathbf{y}^{\mathbf{d}_{\varepsilon_1,\eta}}(F_{\alpha_1}F_{\alpha_2}-\mathbf{y}^{\mathbf{d}_{\alpha_1,\alpha_2}}F_{\alpha_3})+\mathbf{y}^{\mathbf{d}_{\alpha_1,\beta}}(F_{\varepsilon_1}F_{\varepsilon_2}-\mathbf{y}^{\mathbf{d}_{\varepsilon_1,\varepsilon_2}}F_{\varepsilon_3})\\
    &=(\mathbf{y}^{\mathbf{d}_{\varepsilon_1,\eta}}F_{\alpha_1}+\mathbf{y}^{\mathbf{d}_{\alpha_1,\beta}}F_{\varepsilon_1})F_{\alpha_2}-\mathbf{y}^{\mathbf{d}_{\gamma_1,\gamma_2}}(\mathbf{y}^{\mathbf{d}_{\varepsilon_3,\eta}}F_{\alpha_3}+\mathbf{y}^{\mathbf{d}_{\beta,\alpha_3}}F_{\varepsilon_3})\\
    &=F_{\gamma_1}F_{\gamma_2}-\mathbf{y}^{\mathbf{d}_{\gamma_1,\gamma_2}}F_{\gamma_3}.
\end{align*}

Similarly,
\begin{align*}
   \mathbf{g}_{[\gamma]}&=D(\mathbf{g}_{\gamma_1}+\mathbf{g}_{\gamma_2}+\mathbf{e}_n).
\end{align*}

 \end{itemize}
\end{proof}

\begin{remark}
    Applying the formulas derived above for the $F$-polynomial and the $\mathbf{g}$-vector, one observes that the skein relations in $\mathcal{A}_\bullet^\varsigma(\Tilde{T})$ (Proposition \ref{prop: relations in A^sigma}) become identities in $\mathcal{A}_\bullet(T)$.
\end{remark}

\begin{definition}
    Let $[\gamma] \not\subset \Tilde{T}$ be a non-crossing $\varsigma$-orbit of $(\Tilde{\textbf{S}},\Tilde{\textbf{M}})$. We say that $[\gamma]$ is \emph{nice} if either $\operatorname{Res}([\gamma])=\{\gamma_1\}$ or $\operatorname{Res}([\gamma])=\{\gamma_1,\gamma_2\}$, with $\gamma_1$ and $\gamma_2$ crossing at the endpoint $\blacksquare$ in $(\textbf{S},\textbf{M})$.
\end{definition}

\begin{remark}
    If $(\Tilde{\textbf{S}},\Tilde{\textbf{M}})$ is a regular polygon, then every $\varsigma$-orbit is nice.
\end{remark}

\section{Snake graphs and perfect matching Laurent polynomials from surfaces}\label{s2}

In this section, we associate with each nice $\varsigma$-orbit $[\gamma]$ of $(\Tilde{\textbf{S}},\Tilde{\textbf{M}})$ a labeled modified snake graph $\mathcal{G}_{[\gamma]}$ constructed by gluing together the snake graphs corresponding to the arcs of $\operatorname{Res}([\gamma])$. This allows us to obtain the cluster expansion of the cluster variable $x_{[\gamma]}$ of $\mathcal{A}_\bullet^\varsigma(\Tilde{T})$ in terms of perfect matchings of $\mathcal{G}_{[\gamma]}$.

\subsection{Labeled snake graphs from surfaces}\label{sect: sg from poligons} 
We first briefly recall the construction of snake graphs from arcs on triangulated surfaces. We refer the reader to \cite{CSI,MS,MSW2011} for more details.
\begin{definition}[Tile $G_j$]
    Let $T=\{\tau_1, \dots,\tau_n\}$ be a triangulation of a surface $(\textbf{S},\textbf{M})$, and let $\gamma$ be an arc of $(\textbf{S},\textbf{M})$ not contained in $T$. We orient $\gamma$ such that $s$ is its starting point and $t$ its endpoint. Let

\begin{center}
    $p_1,p_2,\dots,p_{d}$,
\end{center}
with $p_j \in \tau_{i_j}$, $j=1,\dots,d$, be the intersections of $\gamma$ with $T$ in order of appearance, and set $p_0=s$ and $p_{d+1}=t$. Let $\Delta_{j-1}$ and $\Delta_j$ be the two triangles of $T$ on each side of $\tau_{i_j}$. The \emph{tile} $G_j$ is the graph with four vertices and five edges, consisting of a square together with one diagonal. Its five edges are in bijection with the five arcs occurring in the triangles $\Delta_{j-1}$ and $\Delta_j$, with the diagonal corresponding to the arc $\tau_{i_j}$. Moreover, this bijection preserves the relative positions of the corresponding arcs, up to the orientation of the two triangles.
\end{definition}

\begin{definition}[Relative orientation]
    Given a planar embedding $\Tilde{G}_j$ of $G_j$, the \emph{relative orientation} $\operatorname{rel}(\Tilde{G}_j,T)$ of $\Tilde{G}_j$ with respect to $T$ is $+1$ (resp. $-1$) if its triangles agree (resp. disagree) in orientation with the corresponding triangles of $T$. 
\end{definition}

 The arcs $\tau_{i_j}$ and $\tau_{i_{j+1}}$ form two edges of the triangle $\Delta_j$ in $T$. We label $\tau_{[j]}$ the third edge of this triangle. 

 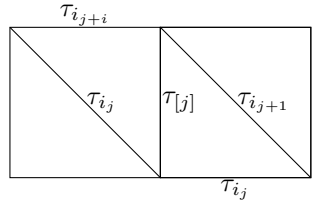
\begin{figure}[h]
    \centering
   \begin{tikzpicture}[scale=2]
  
 
   \draw (-0.5,0) -- (-0.5,1) -- (0.5,1) -- (0.5,0) -- cycle;

\begin{scope}
    \draw (-0.5,0) --node[midway] {} (-0.5,1) --node[midway, above left,xshift=2mm] {} (0.5,1) -- node[midway,right,xshift=-1mm] {}(0.5,0) -- node[midway, below,yshift=1mm] {$\tau_{i_j}$} cycle;
    \draw (-0.5,1) -- node[midway,right,xshift=-1mm] {$\tau_{i_{j+1}}$}(0.5,0);

\end{scope}

\begin{scope}[xshift=-1cm]
    \draw (-0.5,0) -- (-0.5,1) -- node[midway, above,yshift=-1mm] {$\tau_{i_{j+i}}$}(0.5,1) -- node[midway, right,xshift=-1mm] {$\tau_{[j]}$}(0.5,0) --  cycle;
     \draw (-0.5,1) -- node[midway,right,xshift=-1mm] {$\tau_{i_{j}}$}(0.5,0);
\end{scope}

\end{tikzpicture}
\caption{Gluing tiles $\Tilde{G}_j$ and $\Tilde{G}_{j+1}$ along the edge labeled $\tau_{[j]}$.}
\label{gluing}
\end{figure}

\begin{definition}[Snake graph $\mathcal{G}_\gamma$]\label{def_snake_graph}
Let $G_1, \dots , G_d$ be the tiles associated with the intersections of $\gamma$ with $T$. We glue them together in the following way: $\Tilde{G}_{j+1}$ is glued to $\Tilde{G}_j$, along the edge $\tau_{[j]}$, in such a way that $\operatorname{rel}(\Tilde{G}_{j+1},T)\neq \operatorname{rel}(\Tilde{G}_{j},T)$, as in Figure \ref{gluing}. The resulting planar graph is denoted by $\mathcal{G}_\gamma^\Delta$. The \emph{snake graph $\mathcal{G}_\gamma$ associated with $\gamma$} is obtained from $\mathcal{G}_\gamma^\Delta$ by removing the diagonal in each tile. The edges along which we glue two tiles are called \emph{internal}; the other ones are called \emph{external}.   
\end{definition}
\begin{definition}[Perfect matching]
A \emph{perfect matching} of $\mathcal{G}_\gamma$ is a subset $P$ of the edges of $\mathcal{G}_\gamma$ such that each vertex of $\mathcal{G}_\gamma$ is incident to exactly one edge of $P$.    
\end{definition}

\begin{definition}[Minimal and maximal matching of $\mathcal{G}_\gamma$]\label{p-}
    The snake graph $\mathcal{G}_\gamma$ has precisely two perfect matchings which contain only external edges. If $\operatorname{rel}(\Tilde{G}_1,T)=+1$ (resp. $-1$), $e_1$ and $e_2$ are defined to be the two edges of $\mathcal{G}_\gamma^\Delta$ which lie in the counterclockwise (resp. clockwise) direction from the diagonal of $\Tilde{G}_1$. Then $P_-=P_-(\mathcal{G}_\gamma^\Delta)$ is the unique matching which contains only external edges and does not contain edges $e_1$ or $e_2$. $P_-$ is called the \emph{minimal matching}. $P_+=P_+(\mathcal{G}_\gamma^\Delta)$, the \emph{maximal matching}, is the other matching with only external edges.
\end{definition}

Let $P_-\ominus P = (P_- \cup P) \setminus (P_- \cap P)$ be the symmetric difference of the minimal matching $P_-$ and a perfect matching $P$ of $\mathcal{G}_\gamma$. By \cite[Theorem 5.1]{MS}, $P_-\ominus P$ is the set of boundary edges of a subgraph $\mathcal{G}_P$ of $\mathcal{G}_\gamma$, and $\mathcal{G}_P$ is a union of tiles
\begin{center}
    $\mathcal{G}_P=\displaystyle\bigcup_{j \in J} G_j$.
\end{center}
\begin{definition}[Height monomial]\label{def_h(P)_original}
Let $P$ be a perfect matching of $\mathcal{G}_\gamma$. The \emph{height monomial of $P$} is
\begin{center}
    $y(P):=\displaystyle\prod_{j \in J}y_{i_j}$.
\end{center}
Thus $y(P)$ is the product of all $y_{i_j}$ for which the tile $G_j$ lies inside $P_-\ominus P$.
\end{definition}

\begin{notation}
For a tile $G_j$ and a set $E$ of edges of $\mathcal{G}_\gamma$, we write $E|_{G_j}$ for the subset of $E$ consisting of edges of $G_j$.
\end{notation}

\begin{lemma}[{\cite[Lemma 1.7]{ciliberti2}}]\label{lemma_ind_set}
    Let 
    \begin{center}
        $\Tilde{J}=\{ j \mid \text{$(P_- \cup P)|_{G_j}$ contains an external edge of $\mathcal{G}_\gamma$ and $(P_- \cap P)|_{G_j}=\emptyset $}\}$.
    \end{center} Then $\Tilde{J}=J$.
\end{lemma}

\begin{remark}\label{rmk_ind_set}
    It follows from Lemma \ref{lemma_ind_set} that 
    $y(P)$ is the product of all $y_{i_j}$ such that $(P_- \cup P)|_{G_j}$ contains an external edge of $\mathcal{G}_\gamma$ and $(P_- \cap P)|_{G_j}=\emptyset $.
    
\end{remark}

\begin{definition}[Perfect matching polynomial and $\mathbf{g}$-vector of $\mathcal{G}_\gamma$]
    Let $\gamma$ be an arc that is not in $T$, and $\tau_{i_1}, \dots, \tau_{i_d}$ be the arcs of $T$ crossed by $\gamma$. Then the \emph{perfect matching polynomial of} $\mathcal{G}_\gamma$ is
  \begin{center}
      $F_{\mathcal{G}_\gamma}:=\displaystyle\sum_{P}y(P)$,
  \end{center}
  where the sum is over all perfect matchings $P$ of $\mathcal{G}_\gamma$, and the $\mathbf{g}$-$vector$ is
  \begin{center}
      $\mathbf{g}_{\mathcal{G}_\gamma}:=\displaystyle\sum_{e \in P_-(\mathcal{G}_\gamma)}x(e)-\displaystyle\sum_{j=1}^d \mathbf{e}_{i_j}$,
  \end{center}
where, denoting by ${\mathbf{e}_1,\dots,\mathbf{e}_n}$ the standard basis of $\mathbb{Z}^n$, $x(e)=\mathbf{e}_i$ if the label of $e$ is $\tau_i$, and $x(e)=\mathbf{0 \in \mathbb{Z}^n}$ otherwise.
  The definition is extended to any arc by letting $F_{\mathcal{G}_\gamma}:=1$ and $\mathbf{g}_{\mathcal{G}_\gamma}:=\mathbf{e}_i$ if $\gamma=\tau_i \in T$, and $F_{\mathcal{G}_\gamma}:=1$ and $\mathbf{g}_{\mathcal{G}_\gamma}:=\mathbf{0}$ if $\gamma$ is a boundary arc of $(\mathbf{S},\mathbf{M})$. 
\end{definition}

\begin{theorem}[{\cite[Theorem 3.1]{MS}}]\label{f-poly g-vect of gamma}
  Let $T$ be a triangulation of $(\mathbf{S},\mathbf{M})$, and let $\gamma$ be an arc.  Then $F_{\mathcal{G}_\gamma}$ and $\mathbf{g}_{\mathcal{G}_\gamma}$ are the $F$-polynomial $F_\gamma$ and the $\mathbf{g}$-vector $\mathbf{g}_\gamma$, respectively, of the cluster variable $x_\gamma$ of $\mathcal{A}_\bullet(T)$ which corresponds to $\gamma$. 
\end{theorem}

\subsection{Modified snake graphs from $\varsigma$-orbits}

\begin{notation}
    Let $T=\{ \tau_1, \dots, \tau_n \}$ be a triangulation of $(\textbf{S},\textbf{M})$. We denote by $G(i)$ the tile of $\mathcal{G}_{\gamma}$ with diagonal label $\tau_{i}$.
\end{notation}

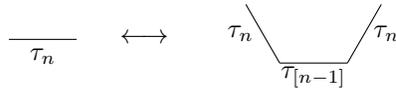
\begin{figure}[h] 
        \centering
         \begin{tikzpicture}[scale=0.9]
\draw (180:1) -- node[midway, left] {$\textcolor{black}{\tau_n}$}(240:1);
\draw (240:1) --node[midway, below,yshift=0.8mm] {$\textcolor{black}{\tau_{[n-1]}}$} (300:1);
\draw (300:1) -- node[midway, right] {$\textcolor{black}{\tau_n}$}(0:1);
\begin{scope}[xshift=-4cm,yshift=0.35cm]
    \draw (240:1) -- node[midway, below] {$\textcolor{black}{\tau_n}$}(300:1);
\end{scope}
\node at (-2.5,-0.5) {$\longleftrightarrow$};
\end{tikzpicture}
       \caption{From $G(n-1)$ to $\hat{G}(n-1)$.}
       \label{fig:modif_sg}
    \end{figure}

\begin{definition}[Labeled modified snake graph $\hat{\mathcal{G}_{\gamma}}$]\label{def_mod_sg_B}
    Let $T=\{ \tau_1, \dots, \tau_n \}$ be a triangulation of $(\textbf{S},\textbf{M})$, such that $\tau_n$ is an arc of a triangle of $T$ whose other two edges are boundary arcs. Let $\gamma$ be an arc of $(\textbf{S},\textbf{M})$ that is not in $T$. We define the \emph{labeled modified snake graph $\hat{\mathcal{G}_{\gamma}}$ associated with $\gamma$} as the labeled snake graph $\mathcal{G}_{\gamma}$ of Definition \ref{def_snake_graph} with the following modifications:
    \begin{itemize}
        \item The edge with label $\tau_n$ in $G(n-1)$ is replaced by three new edges, so that $\hat{G}(n-1)$ has the shape of a hexagon, as in Figure \ref{fig:modif_sg}.
        \item Let $e$ be an internal edge of $\mathcal{G}_\gamma$ contained in the tile $G(n)$, and let $e'$ be the edge of $\hat{G}(n)$ lying on the opposite side of $G(n)$ from $e$. If $l$ is the label of $e$, then $e'$ has the same label $l$.
    \end{itemize}
    
\end{definition}

    \begin{example}
Figure \ref{fig:ex_sg_parts} shows the construction of the labeled modified snake graphs  $\hat{\mathcal{G}}_{\gamma_1}$ and $\hat{\mathcal{G}}_{\gamma_2}$ associated with the arcs $\gamma_1$ and $\gamma_2$ depicted on the right side of Figure \ref{fig:sigma_orbit}. 
\begin{figure}[h] 
    \centering
    \begin{tikzpicture}[scale=0.7]
  
    \begin{scope}[local bounding box=Gamma1]
    
        \draw (-0.5,0) -- (-0.5,1) -- (0.5,1) -- (0.5,0) -- cycle;
        \node at (-3,0.5) {$\mathcal{G}_{\gamma_1}=$};

        \begin{scope}
            \draw (-0.5,0) --node[midway] {} (-0.5,1) --node[midway, above left,xshift=2mm,yshift=-1mm] {\textcolor{blue}{3}} (0.5,1) -- node[midway,right,xshift=-1mm] {\textcolor{blue}{2}}(0.5,0) -- node[midway, below,yshift=1mm] {\textcolor{blue}{5}} cycle;
            \node at (0,0.5) {4};
        \end{scope}
       
        \begin{scope}[xshift=-1cm]
            \draw (-0.5,0) -- (-0.5,1) -- node[midway, above,yshift=-1mm] {\textcolor{blue}{4}}(0.5,1) -- (0.5,0) --  cycle;
            \node at (0,0.5) {5};
        \end{scope}  

        \begin{scope}[xshift=8cm,yshift=0.855cm]
            \node at (-3.5,-0.5) {$\hat{\mathcal{G}}_{\gamma_1}=$};

            \begin{scope}[]
                \filldraw[fill=red!10] (0:1) -- node[midway,right] {\textcolor{blue}{5}} (60:1) -- node[midway,above,yshift=-0.8mm] {\textcolor{blue}{2}} (120:1) -- node[midway,left] {\textcolor{blue}{3}}(180:1) --(240:1)--  node[midway,below] {\textcolor{blue}{5}} (300:1) --node[midway] {\textcolor{blue}{[4]}} (0:1); 
                \node at (0,0) {$4$};
            \end{scope}
       
            \begin{scope}[xshift=-1.5cm,yshift=-0.855cm]
                \filldraw[fill=red!10] (0:1) -- node[pos=0.5, right, xshift=-4mm, yshift=0mm] {} (60:1) -- node[midway] {\textcolor{blue}{4}} (120:1) -- node[midway] {\textcolor{blue}{[4]}}(180:1) -- (0:1);
                \node at (0,0.4) {$5$};
            \end{scope}
        \end{scope}
        
    \end{scope} 

    \begin{scope}[yshift=-5.5cm, local bounding box=Gamma2] 
    
        \draw (-0.5,0) -- (-0.5,1) -- (0.5,1) -- (0.5,0) -- cycle;
        \node at (-3,0.5) {$\mathcal{G}_{\gamma_2}=$};

        \begin{scope}
            \draw (-0.5,0) --node[midway] {} (-0.5,1) --node[midway, above left,xshift=2mm,yshift=-1mm] {\textcolor{blue}{3}} (0.5,1) -- node[midway,right,xshift=-1mm] {\textcolor{blue}{}}(0.5,0) -- node[midway, below,yshift=1mm] {\textcolor{blue}{5}} cycle;
            \node at (0,0.5) {4};
        \end{scope}

        \begin{scope}[xshift=1cm]
            \draw (-0.5,0) -- (-0.5,1) -- node[midway, above,yshift=-1mm] {\textcolor{blue}{1}}(0.5,1) -- (0.5,0) --node[midway, below,yshift=1mm] {\textcolor{blue}{4}}  cycle;
            \node at (0,0.5) {3};
        \end{scope}  
       
        \begin{scope}[xshift=-1cm]
            \draw (-0.5,0) -- (-0.5,1) -- node[midway, above,yshift=-1mm] {\textcolor{blue}{4}}(0.5,1) -- (0.5,0) --  cycle;
            \node at (0,0.5) {5};
        \end{scope}  
        \begin{scope}[xshift=2cm]
            \draw (-0.5,0) -- (-0.5,1) --node[midway, above,yshift=-1mm] {\textcolor{blue}{2}} (0.5,1) -- (0.5,0) -- node[midway, below,yshift=1mm] {\textcolor{blue}{3}} cycle;
            \node at (0,0.5) {1};
        \end{scope} 

        \begin{scope}[xshift=9.5cm, yshift=1.7cm] 
           \filldraw[fill=blue!10] (-1,0) -- node[right,xshift=-1mm] {\textcolor{blue}{4}} (-1,1) -- (-2,1) --node[left,xshift=1mm] {\textcolor{blue}{1}}  (-2,0) --  cycle;
             \node at (-1.5,0.5) {3}; 
        \end{scope}
        \begin{scope}[xshift=9.5cm, yshift=2.7cm]
           \filldraw[fill=blue!10] (-1,0) -- node[right,xshift=-1mm] {\textcolor{blue}{3}} (-1,1) -- (-2,1) --node[left,xshift=1mm] {\textcolor{blue}{2}} (-2,0) -- node[midway,yshift=1mm] {} cycle;
             \node at (-1.5,0.5) {1}; 
        \end{scope}

        \begin{scope}[xshift=5cm,yshift=0.855cm] 
             \node at (-0.5,-0.5) {$\hat{\mathcal{G}}_{\gamma_2}=$};
                
            \begin{scope}[xshift=3cm]
                \filldraw[fill=blue!10] (0:1) -- node[midway,right] {\textcolor{blue}{5}} (60:1) -- node[midway,above] {} (120:1) -- node[midway,left] {\textcolor{blue}{3}}(180:1) --(240:1)--node[midway,below] {\textcolor{blue}{5}} (300:1) --node[midway] {\textcolor{blue}{[4]}} (0:1); 
                \node at (0,0) {$4$};
            \end{scope}
            
            \begin{scope}[xshift=1.5cm,yshift=-0.855cm]
                \filldraw[fill=blue!10] (0:1) -- node[pos=0.5, right, xshift=-4mm, yshift=0mm] {} (60:1) -- node[midway] {\textcolor{blue}{4}} (120:1) -- node[midway] {\textcolor{blue}{[4]}}(180:1) -- (0:1);
                \node at (0,0.4) {$5$};
            \end{scope}

        \end{scope}
    
    \end{scope} 

    \end{tikzpicture}
    \caption{Constructing $\hat{\mathcal{G}}_{\gamma_1}$ and $\hat{\mathcal{G}}_{\gamma_2}$ from $\operatorname{Res}([\gamma])=\{\gamma_1,\gamma_2\}$.}
    \label{fig:ex_sg_parts}
\end{figure}
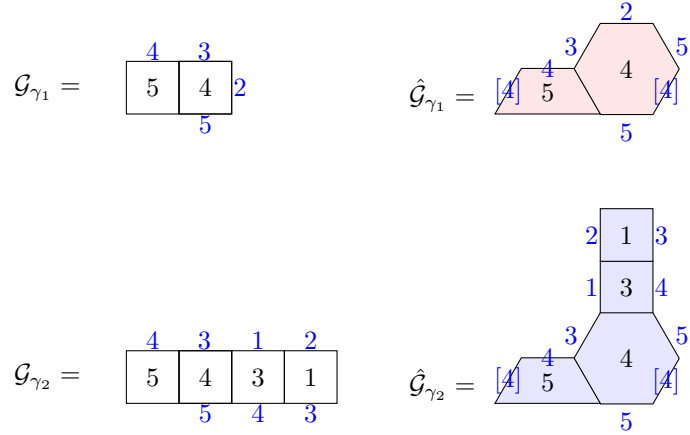 
    \end{example}
\begin{remark}
    In $\hat{\mathcal{G}_{\gamma}}$, unlike $\mathcal{G}_{\gamma}$, $\tau_{[n-1]}$ can also be the label of an external edge. This is the edge along which we will glue the labeled modified snake graphs of some arcs to construct the labeled modified snake graphs associated with $\varsigma$-orbits. See Definition \ref{lab_mod_sg_orbit}.
\end{remark}

\begin{remark}\label{rmk_poset_iso}
By \cite[Section 3, Theorem 2]{propp2024latticestructureorientationsgraphs}, the sets $Match(\mathcal{G}\gamma)$ and $Match(\hat{\mathcal{G}}\gamma)$ of perfect matchings of $\mathcal{G}\gamma$ and $\hat{\mathcal{G}}\gamma$, respectively, are distributive lattices, with partial order given by flips of edges within a tile. We refer the reader to \cite[Section 3, Theorem 2]{propp2024latticestructureorientationsgraphs} for the precise convention adopted. The operation $f:Match(\mathcal{G}_\gamma)\to Match(\hat{\mathcal{G}}_\gamma)$ defined in Figure \ref{fig:poset_iso} is a poset preserving isomorphism.
\end{remark}
        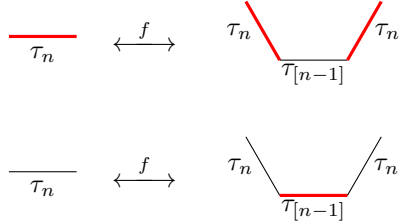
\begin{figure}[h] 
        \centering
         \begin{tikzpicture}[scale=0.9]
\draw[red,very thick] (180:1) -- node[midway, left] {$\textcolor{black}{\tau_n}$}(240:1);
\draw (240:1) --node[midway, below,yshift=0.8mm] {$\textcolor{black}{\tau_{[n-1]}}$} (300:1);
\draw[red,very thick] (300:1) -- node[midway, right] {$\textcolor{black}{\tau_n}$}(0:1);
\begin{scope}[xshift=-4cm,yshift=0.35cm]
    \draw[red,very thick] (240:1) -- node[midway, below] {$\textcolor{black}{\tau_n}$}(300:1);
\end{scope}
\node at (-2.5,-0.5) {$\xlongleftrightarrow{f}$};
\begin{scope}[yshift=-2cm]
     \draw (180:1) -- node[midway, left] {$\textcolor{black}{\tau_n}$}(240:1);
\draw[red,very thick] (240:1) --node[midway, below,yshift=0.8mm] {$\textcolor{black}{\tau_{[n-1]}}$} (300:1);
\draw (300:1) -- node[midway, right] {$\textcolor{black}{\tau_n}$}(0:1);
\begin{scope}[xshift=-4cm,yshift=0.35cm]
    \draw (240:1) -- node[midway, below] {$\textcolor{black}{\tau_n}$}(300:1);
\end{scope}
\node at (-2.5,-0.5) {$\xlongleftrightarrow{f}$};    
\end{scope}
\end{tikzpicture}
       \caption{Poset preserving isomorphism between $Match(\mathcal{G}_\gamma)$ and $Match(\hat{\mathcal{G}}_\gamma)$.}
       \label{fig:poset_iso}
    \end{figure}

\begin{definition}[Labeled modified snake graph $\mathcal{G}_{[\gamma]}$]\label{lab_mod_sg_orbit}
    Let $\Tilde{T} = \{\tau_i\}_{i=1}^{n-1} \sqcup \{\tau_n\} \sqcup \{\tau_i'\}_{i=1}^{n-1}$ be an admissible $\varsigma$-invariant triangulation of a surface $(\Tilde{\textbf{S}},\Tilde{\textbf{M}})$ endowed with an orientation-preserving diffeomorphism $\varsigma$ of order 2. Assume further that $\tau_n$ and $\tau_{n-1}$ are edges of a triangle of $\Tilde{T}$ whose third edge is a boundary arc. Let $[\gamma]$ be a nice $\varsigma$-orbit not contained in $\Tilde{T}$. We define the \emph{labeled modified snake graph $\mathcal{G}_{[\gamma]}$ associated with $[\gamma]$} as follows:
    \begin{itemize}
        \item If $\operatorname{Res}([\gamma])=\{ \gamma_1 \}$, then $\mathcal{G}_{[\gamma]}:=\hat{\mathcal{G}}_{\gamma_1}$;
        \item If $\operatorname{Res}([\gamma])=\{ \gamma_1,\gamma_2\}$, with $\gamma_1$ counterclockwise (resp. clockwise) from $\gamma_2$ at the marked point $\blacksquare$ if $\tau_{n-1}$ is counterclockwise (resp. clockwise) from $\tau_n$, then $\mathcal{G}_{[\gamma]}$ is obtained by gluing $\hat{\mathcal{G}}_{\gamma_2}(n)$ to $\hat{\mathcal{G}}_{\gamma_1}(n-1)$ along $\tau_{[n-1]}$, following the gluing rule recalled in Section \ref{sect: sg from poligons}. If both $\hat{\mathcal{G}}_{\gamma_1}$ and $\hat{\mathcal{G}}_{\gamma_2}$ contain a tile with diagonal label $n-1$, we can assume, after rotating $\mathcal{G}_{[\gamma]}$ if necessary, that $\hat{\mathcal{G}}_{\gamma_2}(n)$ is attached to the right of $\hat{\mathcal{G}}_{\gamma_1}(n-1)$. We then add an edge with label $n-1$ from the top right vertex of $\hat{\mathcal{G}}_{\gamma_1}(n)$ to the top left vertex of $\hat{\mathcal{G}}_{\gamma_2}(n-1)$, as in Figure \ref{additional_edge}.

    \end{itemize}

 The edges of $\mathcal{G}_{[\gamma]}$ along which two tiles are glued are called \emph{internal}, while the remaining edges are called \emph{external}.   
\end{definition}
        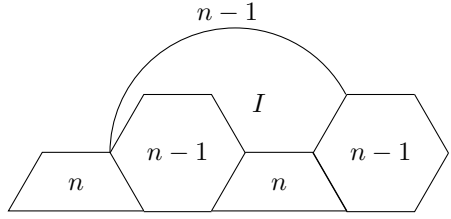
\begin{figure}[h] 
    \centering
   \begin{tikzpicture}[scale=0.9]
   \begin{scope}[xshift=7.8cm, yshift=-2cm]
   \draw (180:1) arc (180:28:1.85) node[midway, above,xshift=3mm] {$n-1$};
   \draw (0:1) -- node[midway,right] {}(60:1) -- (120:1) -- node[midway,left] {}(180:1);
    \draw (300:1) -- node[midway,below] {}(240:1);

\node at (0,0) {$n-1$};

\begin{scope}[xshift=3cm]
\draw[] (0:1) -- node[pos=0.5, right, xshift=-4mm, yshift=0mm] {} (60:1);
    \draw (60:1) -- node[midway,above] {} (120:1) -- (180:1) --(240:1)-- (300:1) -- (0:1); 
    \node at (0,0) {$n-1$};
\end{scope}

\begin{scope}[xshift=-1.5cm,yshift=-0.855cm]
    \draw[] (0:1) -- node[pos=0.5, right, xshift=-4mm, yshift=0mm] {} (60:1);
    \draw (60:1) -- node[midway,above] {} (120:1) -- (180:1) -- (0:1);

    \node at (0,0.4) {$n$};
    
\end{scope}
\begin{scope}[xshift=1.5cm,yshift=-0.855cm]
    \draw (0:1) -- node[midway] {}(60:1) --  (120:1); 
    \draw[] (120:1) -- node[pos=0.5, left,xshift=2mm] {}(180:1);
   \draw (180:1) -- node[midway,below] {}(0:1);
    \node at (-0.3,1.6) {$I$};
    \node at (0,0.4) {$n$};
    
\end{scope}
\end{scope}
   \end{tikzpicture}
    \caption{Additional edge from the top right vertex of $\hat{\mathcal{G}}_{\gamma_1}(n)$ to the top left vertex of $\hat{\mathcal{G}}_{\gamma_2}(n-1)$.}
    \label{additional_edge}
\end{figure}

\begin{example}\label{ex_lab_snake_graph2_B_computation}
Figure \ref{fig:ex_sg_gluing} shows the labeled modified snake graph $\mathcal{G}_{[\gamma]}$ of the $\varsigma$-orbit $[\gamma]$ in Figure \ref{fig:sigma_orbit}, obtained by gluing the graphs $\hat{\mathcal{G}}_{\gamma_1}$ (in red) and $\hat{\mathcal{G}}_{\gamma_2}$ (in blue) constructed in Figure \ref{fig:ex_sg_parts}, and adding an edge with label 4 from the top right vertex of $\hat{\mathcal{G}}_{\gamma_1}(5)$ to the top left vertex of $\hat{\mathcal{G}}_{\gamma_2}(4)$.

\begin{figure}[h] 
\centering
\begin{tikzpicture}[scale=0.9]

   \draw (-1,0) .. controls (-6,1.5) and (0,2.3) .. (2.5,0.85) node[midway, above,xshift=3mm] {\textcolor{blue}{4}};

                \filldraw[fill=red!10] (0:1) -- node[midway,right] {\textcolor{blue}{5}} (60:1) -- node[midway,above,yshift=-0.8mm] {\textcolor{blue}{2}} (120:1) -- node[midway,left] {\textcolor{blue}{3}}(180:1) --(240:1)--  node[midway,below] {\textcolor{blue}{5}} (300:1) --node[midway] {\textcolor{blue}{[4]}} (0:1); 
                \node at (0,0) {$4$};
\node at (-4,0) {$\mathcal{G}_{[\gamma]}=$};
\begin{scope}[xshift=4.5cm, yshift=0.8cm]
 \filldraw[fill=blue!10] (-1,0) -- node[right,xshift=-1mm] {\textcolor{blue}{4}} (-1,1) -- (-2,1) --node[left,xshift=1mm] {\textcolor{blue}{1}}  (-2,0) --  cycle;
             \node at (-1.5,0.5) {3}; 
\end{scope}
\begin{scope}[xshift=4.5cm, yshift=1.8cm]
         \filldraw[fill=blue!10] (-1,0) -- node[right,xshift=-1mm] {\textcolor{blue}{3}} (-1,1) -- (-2,1) --node[left,xshift=1mm] {\textcolor{blue}{2}} (-2,0) -- node[midway,yshift=1mm] {} cycle;
             \node at (-1.5,0.5) {1};  
\end{scope}
\begin{scope}[xshift=3cm]
                \filldraw[fill=blue!10] (0:1) -- node[midway,right] {\textcolor{blue}{5}} (60:1) -- node[midway,above] {} (120:1) -- node[midway,left] {\textcolor{blue}{3}}(180:1) --(240:1)--node[midway,below] {\textcolor{blue}{5}} (300:1) --node[midway] {\textcolor{blue}{[4]}} (0:1); 
                \node at (0,0) {$4$};
\end{scope}

\begin{scope}[xshift=-1.5cm,yshift=-0.855cm]
                \filldraw[fill=red!10] (0:1) -- node[pos=0.5, right, xshift=-4mm, yshift=0mm] {} (60:1) -- node[midway] {\textcolor{blue}{4}} (120:1) -- node[midway] {\textcolor{blue}{[4]}}(180:1) -- (0:1);
                \node at (0,0.4) {$5$};
\end{scope}
\begin{scope}[xshift=1.5cm,yshift=-0.855cm]
                \filldraw[fill=blue!10] (0:1) -- node[pos=0.5, right, xshift=-4mm, yshift=0mm] {} (60:1) -- node[midway] {\textcolor{blue}{4}} (120:1) -- node[midway] {\textcolor{blue}{[4]}}(180:1) -- (0:1);
                \node at (0,0.4) {$5$};
    
\end{scope}
\end{tikzpicture}
\caption{Gluing of $\hat{\mathcal{G}}_{\gamma_1}$ and $\hat{\mathcal{G}}_{\gamma_2}$ that gives $\mathcal{G}_{[\gamma]}$.}
\label{fig:ex_sg_gluing}
  \end{figure}
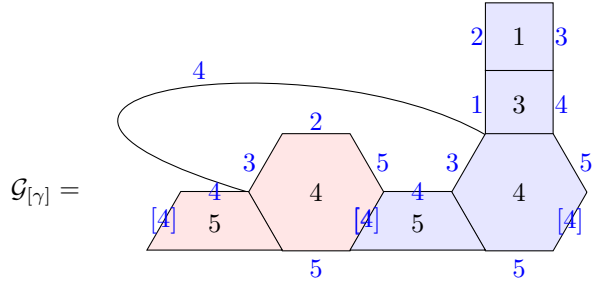
\end{example}
\begin{definition}[Minimal matching of $\mathcal{G}_{[\gamma]}$]\label{def:min matching}
    Let $\mathcal{G}_{[\gamma]}$ be a labeled modified snake graph. We define $P_-(\mathcal{G}_{[\gamma]}) \in Match(\mathcal{G}_{[\gamma]})$ in the following way: 
    \begin{itemize}
        \item if $\mathcal{G}_{[\gamma]}=\hat{\mathcal{G}}_\gamma$, we define $P_-(\mathcal{G}_{[\gamma]}):=f(P_-(\mathcal{G}_\gamma))$, where $f$ is the bijection of Remark \ref{rmk_poset_iso};
        \item if $\mathcal{G}_{[\gamma]}$ is obtained by gluing $\hat{\mathcal{G}}_{\gamma_2}$ to $\hat{\mathcal{G}}_{\gamma_1}$, we define $P_-(\mathcal{G}_{[\gamma]}):=f(P_-(\mathcal{G}_{\gamma_1}))\cup f(P_-(\mathcal{G}_{\gamma_2}))$.
    \end{itemize}
\end{definition}
\begin{example}
    Figure \ref{fig:min pm} shows the minimal perfect matching $P_-(\mathcal{G}_{[\gamma]})$ of the labeled modified snake graph $\mathcal{G}_{[\gamma]}$ in Figure \ref{fig:ex_sg_gluing}. 
    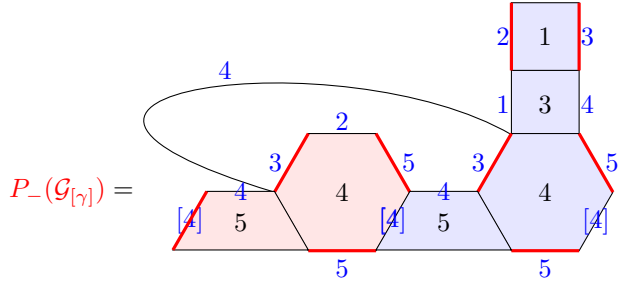
\begin{figure}[h] 
\centering
\begin{tikzpicture}[scale=0.9]

   \draw (-1,0) .. controls (-6,1.5) and (0,2.3) .. (2.5,0.85) node[midway, above,xshift=3mm] {\textcolor{blue}{4}};

                \filldraw[fill=red!10] (0:1) -- node[midway,right] {\textcolor{blue}{5}} (60:1) -- node[midway,above,yshift=-0.8mm] {\textcolor{blue}{2}} (120:1) -- node[midway,left] {\textcolor{blue}{3}}(180:1) --(240:1)--  node[midway,below] {\textcolor{blue}{5}} (300:1) --node[midway] {\textcolor{blue}{[4]}} (0:1); 
                \node at (0,0) {$4$};
                
                \draw[red,very thick] (120:1) -- (180:1);
                \draw[red,very thick] (240:1) -- (300:1);
                \draw[red,very thick] (0:1) -- (60:1);
\node at (-4,0) {$\textcolor{red}{P_-(\mathcal{G}_{[\gamma]})}=$};
\begin{scope}[xshift=4.5cm, yshift=0.8cm]
 \filldraw[fill=blue!10] (-1,0) -- node[right,xshift=-1mm] {\textcolor{blue}{4}} (-1,1) -- (-2,1) --node[left,xshift=1mm] {\textcolor{blue}{1}}  (-2,0) --  cycle;
             \node at (-1.5,0.5) {3}; 
\end{scope}
\begin{scope}[xshift=4.5cm, yshift=1.8cm]
         \filldraw[fill=blue!10] (-1,0) -- node[right,xshift=-1mm] {\textcolor{blue}{3}} (-1,1) -- (-2,1) --node[left,xshift=1mm] {\textcolor{blue}{2}} (-2,0) -- node[midway,yshift=1mm] {} cycle;
             \node at (-1.5,0.5) {1};  
               \draw[red,very thick] (-2,1) -- (-2,0);
                \draw[red,very thick] (-1,0) -- (-1,1);
\end{scope}
\begin{scope}[xshift=3cm]
                \filldraw[fill=blue!10] (0:1) -- node[midway,right] {\textcolor{blue}{5}} (60:1) -- node[midway,above] {} (120:1) -- node[midway,left] {\textcolor{blue}{3}}(180:1) --(240:1)--node[midway,below] {\textcolor{blue}{5}} (300:1) --node[midway] {\textcolor{blue}{[4]}} (0:1); 
                  \draw[red,very thick] (120:1) -- (180:1);
                \draw[red,very thick] (240:1) -- (300:1);
                \draw[red,very thick] (0:1) -- (60:1);
                \node at (0,0) {$4$};
\end{scope}

\begin{scope}[xshift=-1.5cm,yshift=-0.855cm]
                \filldraw[fill=red!10] (0:1) -- node[pos=0.5, right, xshift=-4mm, yshift=0mm] {} (60:1) -- node[midway] {\textcolor{blue}{4}} (120:1) -- node[midway] {\textcolor{blue}{[4]}}(180:1) -- (0:1);
                \node at (0,0.4) {$5$};
                \draw[red,very thick] (120:1) -- (180:1);
\end{scope}
\begin{scope}[xshift=1.5cm,yshift=-0.855cm]
                \filldraw[fill=blue!10] (0:1) -- node[pos=0.5, right, xshift=-4mm, yshift=0mm] {} (60:1) -- node[midway] {\textcolor{blue}{4}} (120:1) -- node[midway] {\textcolor{blue}{[4]}}(180:1) -- (0:1);
                \node at (0,0.4) {$5$};
    
\end{scope}
\end{tikzpicture}
\caption{Minimal perfect matching of $\mathcal{G}_{[\gamma]}$.}
\label{fig:min pm}
  \end{figure}
\end{example}
We extend the definition of height monomial $y(P)$ of a perfect matching $P$ of $\mathcal{G}_{[\gamma]}$ using Remark \ref{rmk_ind_set} to include perfect matchings containing the additional edge. The following definitions are given in \cite{ciliberti2} in the case of regular polygons.
\begin{definition}\label{def_h(P)}
Let $P_-=P_-(\mathcal{G}_{[\gamma]})$, and let $P$ be a perfect matching of $\mathcal{G}_{[\gamma]}$. The \emph{height monomial} of $P$ is
\begin{center}
    $y(P):=\displaystyle\prod_{j \in J}y_{i_j}$,
\end{center}
where the product is over all $i_j$ for which $(P_- \cup P)|_{G_j}$ contains an external edge of $\mathcal{G}_{[\gamma]}$ and $(P_- \cap P)|_{G_j}$ does not contain any edge of $G_j$ with label different from $\tau_n$.
\end{definition}

\begin{definition}[Perfect matching polynomial and $\mathbf{g}$-vector of $\mathcal{G}_{[\gamma]}$]\label{def_f_poly_lab_mod_sg}
    Let $[\gamma]$ be a nice $\varsigma$-orbit not contained in $\Tilde{T}$, and $\tau_{i_1}, \dots, \tau_{i_d}$ be the sequence of  arcs of $T=\operatorname{Res}(\Tilde{T})$ crossed by the  arcs of $\operatorname{Res}([\gamma])$. Then the \emph{perfect matching polynomial of $\mathcal{G}_{[\gamma]}$} is
  \begin{center}
      $F_{\mathcal{G}_{[\gamma]}}:=\displaystyle\sum_{P}y(P)$,
  \end{center}
  where the sum is over all perfect matchings $P$ of $\mathcal{G}_{[\gamma]}$, and the $\mathbf{g}$-$vector$ is
  \begin{center}
      $\mathbf{g}_{\mathcal{G}_{[\gamma]}}:=\displaystyle\sum_{e \in P_-(\mathcal{G}_{[\gamma]})}x(e)-\displaystyle\sum_{j=1}^d \mathbf{e}_{i_j}$,
  \end{center}
where, denoting by ${\mathbf{e}_1,\dots,\mathbf{e}_n}$ the standard basis of $\mathbb{Z}^n$, $x(e)=\mathbf{e}_i$ if the label of $e$ is $\tau_i$, and $x(e)=\mathbf{0 \in \mathbb{Z}^n}$ otherwise. The definition is extended to any $\varsigma$-orbit by letting $F_{\mathcal{G}_{[\gamma]}}:=1$ and $\mathbf{g}_{\mathcal{G}_{[\gamma]}}:=\mathbf{e}_i$ if $[\gamma]=[\tau_i] \in \Tilde{T}$, and $F_{\mathcal{G}_{[\gamma]}}:=1$ and $\mathbf{g}_{\mathcal{G}_{[\gamma]}}:=\mathbf{0}$ if $\gamma$ is a boundary arc of $(\Tilde{\textbf{S}},\Tilde{\textbf{M}})$. 
\end{definition}
\begin{remark}
    In the definition of $\mathbf{g}_{\mathcal{G}_{[\gamma]}}$ external labels of $\mathcal{G}_{[\gamma]}$ of the form $\tau_{[i]}$ are ignored.
\end{remark}
\begin{lemma}\label{lemma_bijection}
    Let $T=\{ \tau_1, \dots, \tau_n \}$ be a triangulation of $(\textbf{S},\textbf{M})$, such that $\tau_n$ is an arc of a triangle of $T$ whose other two edges are boundary arcs. Let $\gamma$ be  an arc of $(\textbf{S},\textbf{M})$ that is not in $T$. Then $F_{\hat{\mathcal{G}}_\gamma}=F_{\mathcal{G}_\gamma}$.
\end{lemma}
\begin{proof}
    Consider the bijection $f:Match(\mathcal{G}_\gamma)\to Match(\hat{\mathcal{G}}_\gamma)$ of Remark \ref{rmk_poset_iso}. We have $y(P)=y(f(P))$, for any $P \in Match(\mathcal{G}_\gamma)$. Therefore, $F_{\hat{\mathcal{G}}_\gamma}=F_{\mathcal{G}_\gamma}$. See \cite[Example 3.14]{ciliberti2} for an illustrating example of this proof in the case of polygons.
\end{proof}

\begin{theorem}\label{thm:sg} 
  Let $\Tilde{T} = \{\tau_i\}_{i=1}^{n-1} \sqcup \{\tau_n\} \sqcup \{\tau_i'\}_{i=1}^{n-1}$ be an admissible $\varsigma$-invariant triangulation of a surface $(\Tilde{\textbf{S}},\Tilde{\textbf{M}})$ endowed with an orientation-preserving diffeomorphism $\varsigma$ of order 2. Assume further that $\tau_n$ and $\tau_{n-1}$ are edges of a triangle of $\Tilde{T}$ whose third edge is a boundary arc. Then, for any nice $\varsigma$-orbit $[\gamma]$ of $(\Tilde{\textbf{S}},\Tilde{\textbf{M}})$, $F_{\mathcal{G}_{[\gamma]}}=F_{[\gamma]}$ and $\mathbf{g}_{\mathcal{G}_{[\gamma]}}=\mathbf{g}_{[\gamma]}$.  
\end{theorem}
\begin{example}\label{ex_lab_snake_graph2_B}
We compute the perfect matching polynomial $F_{\mathcal{G}_{[\gamma]}}$ and the $\mathbf{g}$-vector $\mathbf{g}_{\mathcal{G}_{[\gamma]}}$ of the labeled modified snake graph $\mathcal{G}_{[\gamma]}$ of Figure \ref{fig:ex_sg_gluing} (its minimal perfect matching is displayed in Figure \ref{fig:min pm}):

\begin{align*}
    F_{\mathcal{G}_{[\gamma]}}&=y_1y_3y_4^2y_5^2 + 2y_1y_3y_4^2y_5 + y_1y_4^2y_5^2 + y_1y_3y_4^2 + 2y_1y_4^2y_5 + y_4^2y_5^2 + y_1y_3y_4 \\&+ y_1y_4^2 + 2y_1y_4y_5 + 2y_4^2y_5 + 2y_1y_4 + y_4^2 + 2y_4y_5 + y_1 + 2y_4 + 1;
\end{align*}

\begin{align*}
    \mathbf{g}_{\mathcal{G}_{[\gamma]}}=\begin{pmatrix} 0\\ 1\\ 3\\0\\4
    \end{pmatrix}-\begin{pmatrix} 1\\ 0\\ 1\\2\\2
    \end{pmatrix}=\begin{pmatrix} -1\\1\\2\\-2\\2
    \end{pmatrix}.
\end{align*}
By Theorem \ref{thm:sg}, $F_{\mathcal{G}_{[\gamma]}}$ and $\mathbf{g}_{\mathcal{G}_{[\gamma]}}$ are the $F$-polynomial and the $\mathbf{g}$-vector, respectively, of the cluster variable $x_{[\gamma]}$ corresponding to the $\varsigma$-orbit $[\gamma]$ depicted on the left-hand side of Figure \ref{fig:sigma_orbit}, in the skew-symmetrizable cluster algebra with principal coefficients in the $\varsigma$-invariant triangulation of Figure \ref{fig:ex_adm_triang}. 
\end{example}

\begin{remark}
    If $(\Tilde{\textbf{S}},\Tilde{\textbf{M}})$ is a regular polygon, we recover \cite[Theorem 3.18]{ciliberti2}.
\end{remark}

\begin{proof}[Proof of Theorem \ref{thm:sg}]
The proof is based on Theorem \ref{thm:formula cv}. We first prove that $F_{\mathcal{G}_{[\gamma]}}=F_{[\gamma]}$. Let $T=\operatorname{Res}(\Tilde{T})=\{\tau_1, \dots, \tau_n\}$. If $\operatorname{Res}([\gamma])=\{\gamma_1\}$, the statement holds since $F_{\mathcal{G}_{[\gamma]}}=F_{\hat{\mathcal{G}}_{\gamma_1}}=F_{\mathcal{G}_{\gamma_1}}$ (see Lemma \ref{lemma_bijection}). It follows from Theorem \ref{thm:formula cv} and Theorem \ref{f-poly g-vect of gamma} that $F_{\mathcal{G}_{[\gamma]}}=F_{[\gamma]}$. Otherwise, $\operatorname{Res}([\gamma])=\{\gamma_1,\gamma_2\}$. Since $[\gamma]$ is nice, and $\tau_n$ and $\tau_{n-1}$ are edges of a triangle of $\Tilde{T}$ whose third edge is a boundary arc, there are only two cases to consider:

    \begin{itemize}
        \item [(i)] One of $\gamma_1$ and $\gamma_2$, say $\gamma_2$, intersects only $\tau_n$. So $\mathcal{G}_{[\gamma]}$ is of the form displayed in Figure \ref{fig:sg_only_tau_n}. We have that

\begin{align*}
F_{\mathcal{G}_{[\gamma]}}=F_{\hat{\mathcal{G}}_{\gamma_1}}F_{\hat{\mathcal{G}}_{\gamma_2}}-R=F_{\mathcal{G}_{\gamma_1}}F_{\mathcal{G}_{\gamma_2}}-R,
\end{align*}
where $R$ is the sum of the monomials that correspond to unions of perfect matchings of $\hat{\mathcal{G}}_{\gamma_1}$ and perfect matchings of $\hat{\mathcal{G}}_{\gamma_2}$ which do not give perfect matchings of $\mathcal{G}_{[\gamma]}$. They are all of the form displayed in Figure \ref{forbidden_2}.
        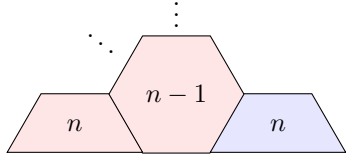
\begin{figure}[h] 
    \centering
   \begin{tikzpicture}[scale=0.9]
   \begin{scope}[xshift=7.8cm, yshift=-2cm]
   \filldraw[fill=red!10] (0:1) -- node[midway,right] {}(60:1) --node[midway, above] {$\vdots$} (120:1) -- node[midway, above left] {$\ddots$}(180:1) -- (240:1) -- (300:1);

\node at (0,0) {$n-1$};

\begin{scope}[xshift=-1.5cm,yshift=-0.855cm]
    \filldraw[fill=red!10] (0:1) -- node[pos=0.5, right, xshift=-4mm, yshift=0mm] {} (60:1) -- node[midway,above] {} (120:1) -- (180:1) -- (0:1);

    \node at (0,0.4) {$n$};
    
\end{scope}
\begin{scope}[xshift=1.5cm,yshift=-0.855cm]
    \filldraw[fill=blue!10] (0:1) -- node[midway] {}(60:1) --  (120:1)-- node[pos=0.5, left,xshift=2mm] {}(180:1)-- node[midway,below] {}(0:1);
    
    \node at (0,0.4) {$n$};
    
\end{scope}
\end{scope}
   \end{tikzpicture}
    \caption{The shape of the graph $\mathcal{G}_{[\gamma]}$ if $\operatorname{Res}([\gamma])=\{\gamma_1,\gamma_2\}$, and $\gamma_2$ intersects only $\tau_n$.}
    \label{fig:sg_only_tau_n}
\end{figure}  

             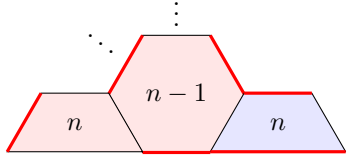
\begin{figure}[h] 
    \centering
   \begin{tikzpicture}[scale=0.9]
   \begin{scope}[xshift=7.8cm, yshift=-2cm]
   \filldraw[fill=red!10] (0:1) -- node[midway,right] {}(60:1) --node[midway, above] {$\vdots$} (120:1) -- node[midway, above left] {$\ddots$}(180:1) -- (240:1) -- (300:1);
\draw[red,very thick] (0:1) -- node[midway] {}(60:1);
\draw[red,very thick] (120:1) -- node[midway] {}(180:1);
\draw[red,very thick] (240:1) -- node[midway] {}(300:1);
\node at (0,0) {$n-1$};

\begin{scope}[xshift=-1.5cm,yshift=-0.855cm]
    \filldraw[fill=red!10] (0:1) -- node[pos=0.5, right, xshift=-4mm, yshift=0mm] {} (60:1)-- node[midway,above] {} (120:1) -- (180:1) -- (0:1);
\draw[red,very thick] (120:1) -- node[midway] {}(180:1);
    
    \node at (0,0.4) {$n$};
    
\end{scope}
\begin{scope}[xshift=1.5cm,yshift=-0.855cm]
    \filldraw[fill=blue!10] (0:1) -- node[midway] {}(60:1) --  (120:1) -- node[pos=0.5, left,xshift=2mm] {}(180:1);
   \draw (180:1) -- node[midway,below] {}(0:1);
    \draw[red,very thick] (120:1) -- node[midway] {}(60:1);
    \draw[red,very thick] (180:1) -- node[midway] {}(0:1);
    \node at (0,0.4) {$n$};
    
\end{scope}
\end{scope}
   \end{tikzpicture}
    \caption{Unions of perfect matchings of $\hat{\mathcal{G}}_{\gamma_1}$ and perfect matchings of $\hat{\mathcal{G}}_{\gamma_2}$ that fail to produce perfect matchings of $\mathcal{G}_{[\gamma]}$ (case where $\gamma_2$ intersects only $\tau_n$).}
    \label{forbidden_2}
\end{figure}

Let $\tilde\gamma_1$ be the  arc of the collapsed surface $(\textbf{S},\textbf{M})$ which intersects the same  arcs of $T$ as $\gamma_1$ but $\tau_n$. To describe $R$, we consider the skein relation in $(\textbf{S},\textbf{M})$ corresponding to the crossing of $\tau_{n-1}$ and $\tilde\gamma_1$. We have two cases to consider:
\begin{itemize}
    \item [1)] The arc $\tau_{[n-1]}$, along which we glue $\hat{\mathcal{G}}_{\gamma_1}$ and $\hat{\mathcal{G}}_{\gamma_2}$, is not in the minimal perfect matching of $\hat{\mathcal{G}}_{\gamma_1}$. So the red edges of $\hat{G}(n-1)$ in Figure \ref{forbidden_2} are in the minimal perfect matching of $\hat{\mathcal{G}}_{\gamma_1}$. It follows that $y_{n-1}$ is a summand of $F_{\mathcal{G}_{\gamma_1}}$. Therefore, $\tau_{n-1}$ must necessarily be counterclockwise from $\tau_n$.
\begin{figure}[h]
    \centering
    \includegraphics[width=0.7\linewidth]{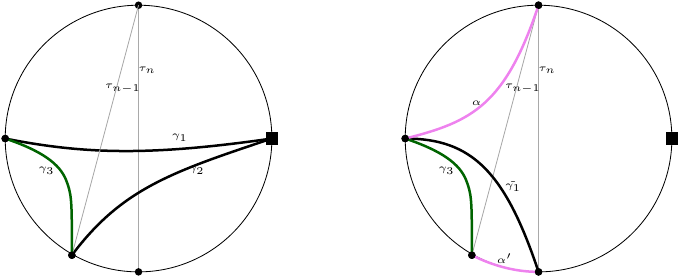}
    \caption{On the left, $\operatorname{Res}([\gamma])=\{\gamma_1,\gamma_2\}$ and the arc $\gamma_3$ resolving the crossing of $\gamma_1$ and $\gamma_2$ at $\blacksquare$; on the right, the resolution of the crossing of $\Tilde{\gamma_1}$ and $\tau_{n-1}$ (case $\gamma_2$ crosses only $\tau_n$, and $\tau_{n-1}$ is counterclockwise from $\tau_n$).}
    \label{fig:proof2_fig1}
\end{figure}
In the notation of Figure \ref{fig:proof2_fig1}, we have that
\begin{equation}\label{eq_3_only_tau_n}
    F_{\tilde\gamma_1}=\mathbf{y}^{\mathbf{d}_{\gamma_3,\tau_n}}F_{\alpha}F_{\alpha'}+F_{\gamma_3}F_{\tau_n}=\mathbf{y}^{\mathbf{d}_{\gamma_3,\tau_n}}F_{\alpha}+F_{\gamma_3}.
\end{equation}
Since the red edges of $\hat{G}(n-1)$ in Figure \ref{forbidden_2} are in the minimal perfect matching of $\hat{\mathcal{G}}_{\gamma_1}$, and so of $\hat{\mathcal{G}}_{\tilde\gamma_1}$, the sum of the monomials which correspond to the perfect matchings of $\hat{\mathcal{G}}_{\tilde\gamma_1}$ which contain those edges in the right hand side of \eqref{eq_3_only_tau_n} is $F_{\gamma_3}$. Thus,
\begin{align*}
    R=y_nF_{\gamma_3}=\mathbf{y}^{\mathbf{d}_{\gamma_1,\gamma_2}}F_{\gamma_3}.
\end{align*}
We observe that, since $\tau_n$ and $\tau_{n-1}$ are sides of a triangle of $T$ whose third edge is a boundary arc, the arc $\gamma_3$ is isotopic to the arc that resolves the crossing of $\gamma_1$ and $\gamma_2$ at $\blacksquare$ (on the left-hand side of Figure \ref{fig:proof2_fig1}). Therefore, by Theorem \ref{thm:formula cv}, $F_{\mathcal{G}_{[\gamma]}}=F_{[\gamma]}$.
\item [2)] The arc $\tau_{[n-1]}$ is in the minimal perfect matching of $\hat{\mathcal{G}}_{\gamma_1}$. The argument in this situation is analogous to case 1).
\end{itemize}
\item [(ii)] Both arcs $\gamma_1$ and $\gamma_2$ intersect $\tau_{n-1}$. It follows that $\mathcal{G}_{[\gamma]}$ is of the form displayed in Figure \ref{fig:sg_both_tau_n-1}. We have that
 \begin{figure}[h] 
    \centering
   \begin{tikzpicture}[scale=0.9]
   \begin{scope}[xshift=7.8cm, yshift=-2cm]
   \draw (180:1) arc (180:28:1.85);

    \filldraw[fill=red!10] (0:1) -- node[midway,right] {}(60:1) --node[midway, above] {$\vdots$} (120:1) -- node[midway, above left] {$\ddots$}(180:1) -- (240:1) -- (300:1);

\node at (0,0) {$n-1$};

\begin{scope}[xshift=3cm]
\filldraw[fill=blue!10] (0:1) -- node[pos=0.5, right, xshift=-4mm, yshift=0mm] {} (60:1) -- node[midway,above] {} (120:1) -- (180:1) --(240:1)-- (300:1) -- (0:1); 
     \draw (0:1) -- node[midway,right] {}(60:1) --node[midway, above] {$\vdots$} (120:1) -- node[midway, above left] {$\ddots$}(180:1);
    \node at (0,0) {$n-1$};
\end{scope}

\begin{scope}[xshift=-1.5cm,yshift=-0.855cm]
    \filldraw[fill=red!10] (0:1) -- node[pos=0.5, right, xshift=-4mm, yshift=0mm] {} (60:1) -- node[midway,above] {} (120:1) -- (180:1) -- (0:1);

    \node at (0,0.4) {$n$};
    
\end{scope}
\begin{scope}[xshift=1.5cm,yshift=-0.855cm]
    \filldraw[fill=blue!10] (0:1) -- node[midway] {}(60:1) --  (120:1) -- node[pos=0.5, left,xshift=2mm] {}(180:1) -- node[midway,below] {}(0:1);
    
    \node at (0,0.4) {$n$};
    
\end{scope}
\end{scope}
   \end{tikzpicture}
\caption{The shape of the graph $\mathcal{G}_{[\gamma]}$ if $\operatorname{Res}([\gamma])=\{\gamma_1,\gamma_2\}$, and both $\gamma_1$ and $\gamma_2$ intersect $\tau_{n-1}$.}
    \label{fig:sg_both_tau_n-1}
\end{figure}
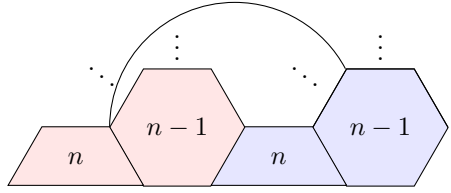
\begin{equation*}
    F_{\mathcal{G}_{[\gamma]}}=F_{\hat{\mathcal{G}}_{\gamma_1}}F_{\hat{\mathcal{G}}_{\gamma_2}}-R + S=F_{\mathcal{G}_{\gamma_1}}F_{\mathcal{G}_{\gamma_2}}-R+S,
\end{equation*}
where $R$ is the sum of the monomials which correspond to unions of perfect matchings of $\hat{\mathcal{G}}_{\gamma_1}$ and perfect matchings of $\hat{\mathcal{G}}_{\gamma_2}$ which are not perfect matchings of $\mathcal{G}_{[\gamma]}$ (all of the form displayed in Figure \ref{forbidden_3}); while $S$ is the sum of the monomials which correspond to perfect matchings of $\mathcal{G}_{[\gamma]}$ which contain the additional edge (all of the form displayed in Figure \ref{forbidden_4}).

 \begin{figure}[h] 
    \centering
   \begin{tikzpicture}[scale=0.9]
   \begin{scope}[xshift=7.8cm, yshift=-2cm]
   \draw (180:1) arc (180:28:1.85);
\node at (-0.5,-1.2) {$\hat{\mathcal{G}}_{\gamma_1}$};
\node at (2.5,-1.2) {$\hat{\mathcal{G}}_{\gamma_2}$};
    \filldraw[fill=red!10] (0:1) -- node[midway,right] {}(60:1) --node[midway, above] {$\vdots$} (120:1) -- node[midway, above left] {$\ddots$}(180:1) -- (240:1) -- (300:1);
\draw[red,very thick] (0:1) -- node[midway] {}(60:1);
\draw[red,very thick] (120:1) -- node[midway] {}(180:1);
\draw[red,very thick] (240:1) -- node[midway] {}(300:1);
\node at (0,0) {$n-1$};

\begin{scope}[xshift=3cm]
\filldraw[fill=blue!10] (0:1) -- node[pos=0.5, right, xshift=-4mm, yshift=0mm] {} (60:1)-- node[midway,above] {} (120:1) -- (180:1) --(240:1)-- (300:1) -- (0:1); 
     \draw (0:1) -- node[midway,right] {}(60:1) --node[midway, above] {$\vdots$} (120:1) -- node[midway, above left] {$\ddots$}(180:1);
    \node at (0,0) {$n-1$};
    \draw[red,very thick] (120:1) -- node[midway] {}(60:1);
\draw[red,very thick] (0:1) -- node[midway] {}(300:1);
\end{scope}

\begin{scope}[xshift=-1.5cm,yshift=-0.855cm]
    \filldraw[fill=red!10] (0:1) -- node[pos=0.5, right, xshift=-4mm, yshift=0mm] {} (60:1) -- node[midway,above] {} (120:1) -- (180:1) -- (0:1);

 \draw[red,very thick] (120:1) -- node[midway] {}(180:1);   
    \node at (0,0.4) {$n$};
    
\end{scope}
\begin{scope}[xshift=1.5cm,yshift=-0.855cm]
    \filldraw[fill=blue!10] (0:1) -- node[midway] {}(60:1) --  (120:1) -- node[pos=0.5, left,xshift=2mm] {}(180:1) -- node[midway,below] {}(0:1);
  \draw[red,very thick] (0:1) -- node[midway] {}(180:1);
   \draw[red,very thick] (120:1) -- node[midway] {}(60:1);
    \node at (0,0.4) {$n$};
    
\end{scope}
\end{scope}
   \end{tikzpicture}
    \caption{Unions of perfect matchings of $\hat{\mathcal{G}}_{\gamma_1}$ and perfect matchings of $\hat{\mathcal{G}}_{\gamma_2}$ that fail to produce perfect matchings of $\mathcal{G}_{[\gamma]}$ (case where both $\gamma_1$ and $\gamma_2$ intersect $\tau_{n-1}$).}
    \label{forbidden_3}
\end{figure}
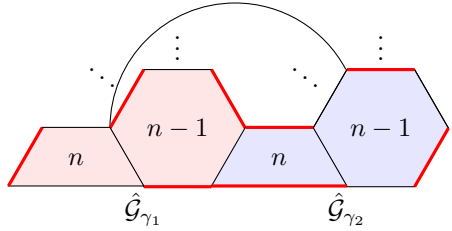

 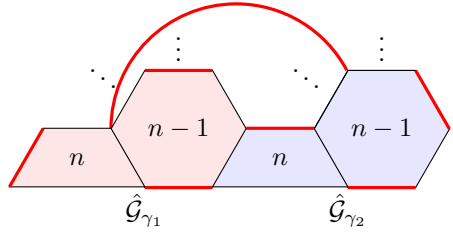
\begin{figure}[h] 
    \centering
   \begin{tikzpicture}[scale=0.9]
   \begin{scope}[xshift=7.8cm, yshift=-2cm]
   \draw[red,very thick] (180:1) arc (180:28:1.85);

    \filldraw[fill=red!10] (0:1) -- node[midway,right] {}(60:1) --node[midway, above] {$\vdots$} (120:1) -- node[midway, above left] {$\ddots$}(180:1)-- (240:1) -- (300:1);
\draw[red,very thick] (120:1) -- node[midway] {}(60:1);
\draw[red,very thick] (240:1) -- node[midway] {}(300:1);
\node at (0,0) {$n-1$};
\node at (-0.5,-1.2) {$\hat{\mathcal{G}}_{\gamma_1}$};
\node at (2.5,-1.2) {$\hat{\mathcal{G}}_{\gamma_2}$};
\begin{scope}[xshift=3cm]
\filldraw[fill=blue!10] (0:1) -- node[pos=0.5, right, xshift=-4mm, yshift=0mm] {} (60:1) -- node[midway,above] {} (120:1) -- (180:1) --(240:1)-- (300:1) -- (0:1); 
     \draw (0:1) -- node[midway,right] {}(60:1) --node[midway, above] {$\vdots$} (120:1) -- node[midway, above left] {$\ddots$}(180:1);
    \node at (0,0) {$n-1$};
    \draw[red,very thick] (0:1) -- node[midway] {}(60:1);
\draw[red,very thick] (240:1) -- node[midway] {}(300:1);
    
\end{scope}

\begin{scope}[xshift=-1.5cm,yshift=-0.855cm]
    \filldraw[fill=red!10] (0:1) -- node[pos=0.5, right, xshift=-4mm, yshift=0mm] {} (60:1)-- node[midway,above] {} (120:1) -- (180:1) -- (0:1);

 \draw[red,very thick] (120:1) -- node[midway] {}(180:1);   
    \node at (0,0.4) {$n$};
    
\end{scope}
\begin{scope}[xshift=1.5cm,yshift=-0.855cm]
    \filldraw[fill=blue!10] (0:1) -- node[midway] {}(60:1) --  (120:1)-- node[pos=0.5, left,xshift=2mm] {}(180:1) -- node[midway,below] {}(0:1);
   \draw[red,very thick] (120:1) -- node[midway] {}(60:1);
    \node at (0,0.4) {$n$};
    
\end{scope}
\end{scope}
   \end{tikzpicture}
    \caption{Perfect matchings of $\mathcal{G}_{[\gamma]}$ which contain the additional edge from the top right vertex of $\hat{\mathcal{G}}_{\gamma_1}(n)$ to the top left vertex of $\hat{\mathcal{G}}_{\gamma_2}(n-1)$.}
    \label{forbidden_4}
\end{figure}
Let $\tilde\gamma_1$ (resp. $\tilde\gamma_2$) be the  arc of $(\textbf{S},\textbf{M})$ that intersects the same arcs of $T$ as $\gamma_1$ (resp. $\gamma_2$) but $\tau_n$. To determine $R$ and $S$, we consider the skein relations corresponding to the crossings of $\tau_{n-1}$ and $\tilde\gamma_1$, and of $\tau_{n-1}$ and $\tilde\gamma_2$. We have two cases to consider.
\begin{figure}[h]
    \centering
    \includegraphics[width=0.7\linewidth]{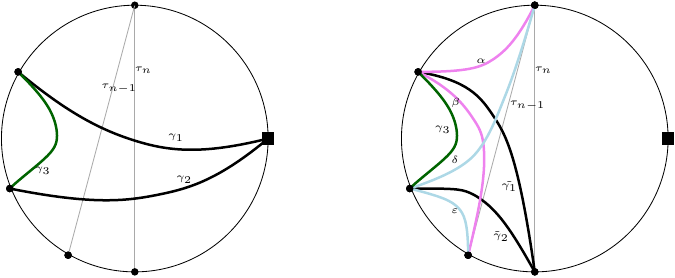}
\caption{On the left, $\operatorname{Res}([\gamma])=\{\gamma_1,\gamma_2\}$ and the arc $\gamma_3$ resolving the crossing of $\gamma_1$ and $\gamma_2$ at $\blacksquare$; on the right, the resolution of the crossing of the arcs $\beta$ and $\delta$ obtained in the resolution of the crossing of $\Tilde{\gamma}_1$ and $\tau_{n-1}$, and of $\Tilde{\gamma_2}$ and $\tau_{n-1}$, respectively (case both $\gamma_1$ and $\gamma_2$ cross $\tau_{n-1}$, and $\tau_{n-1}$ is counterclockwise from $\tau_n$).}
    \label{fig:proof2_fig2}
\end{figure}
\begin{itemize}
    \item [1)] The arc $\tau_{[n-1]}$ is not in the minimal perfect matching of $\hat{\mathcal{G}}_{\gamma_1}$. So the red edges of $\hat{\mathcal{G}}_{\gamma_1}(n-1)$ in Figure \ref{forbidden_3} are in the minimal perfect matching of $\hat{\mathcal{G}}_{\gamma_1}$. It follows that $y_{n-1}$ is a summand of $F_{\mathcal{G}_{\gamma_1}}$. Therefore, $\tau_{n-1}$ must necessarily be counterclockwise from $\tau_n$. In the notation of Figure \ref{fig:proof2_fig2}, we have 
\begin{equation}\label{eq_3}
    F_{\tilde\gamma_1}=\mathbf{y}^{\mathbf{d}_{\beta,\tau_n}}F_{\alpha}+F_{\beta},
\end{equation}

and
\begin{equation}\label{eq_5}
    F_{\tilde\gamma_2}=\mathbf{y}^{\mathbf{d}_{\varepsilon,\tau_n}}F_{\delta}+F_{\varepsilon}.
\end{equation}

Since the red edges of $\hat{\mathcal{G}}_{\gamma_1}(n-1)$ in Figure \ref{forbidden_3} are in the minimal perfect matching of $\hat{\mathcal{G}}_{\gamma_1}$, and so of $\hat{\mathcal{G}}_{\tilde\gamma_1}$, the sum of the monomials corresponding to the perfect matchings of $\hat{\mathcal{G}}_{\tilde\gamma_1}$ that contain those edges on the right hand side of \eqref{eq_3} is $F_\beta$. Since the red edges of $\hat{\mathcal{G}}_{\gamma_1}(n-1)$ in Figure \ref{forbidden_3} are in the minimal perfect matching of $\hat{\mathcal{G}}_{\gamma_1}$, it also follows that the red edges of $\hat{\mathcal{G}}_{\gamma_2}(n-1)$ in Figure \ref{forbidden_3} are not in the minimal perfect matching of $\hat{\mathcal{G}}_{\gamma_2}$. So the sum of the monomials corresponding to the perfect matchings of $\hat{\mathcal{G}}_{\tilde\gamma_2}$ that contain those edges on the right hand side of \eqref{eq_5} is $\mathbf{y}^{\mathbf{d}_{\varepsilon,\tau_n}}F_{\delta}$. Thus,
\begin{align*}
    R=y_n\mathbf{y}^{\mathbf{d}_{\varepsilon,\tau_n}}F_{\delta}F_\beta.
\end{align*}
Similarly,
\begin{align*}
S=y_n\mathbf{y}^{\mathbf{d}_{\beta,\tau_n}}F_{\alpha}F_{\varepsilon}.
\end{align*}

Finally, we consider the skein relation corresponding to the crossing of $\beta$ and $\delta$. We have that
\begin{align*}
F_{\beta}F_{\delta}=\mathbf{y}^{\mathbf{d}_{\alpha,\varepsilon}}F_{\gamma_3}+ \mathbf{y}^{\mathbf{d}_{\gamma_3,\tau_{n-1}}}F_{\alpha}F_{\varepsilon}.
\end{align*}

Therefore,
\begin{align*}
    -R+S=-\mathbf{y}^{\mathbf{d}_{\gamma_1,\gamma_2}}F_{\gamma_3}.
\end{align*}

We observe that, since $\tau_n$ and $\tau_{n-1}$ are sides of a triangle of $T$ whose third edge is a boundary arc, the arc $\gamma_3$ is isotopic to the arc that resolves the crossing of $\gamma_1$ and $\gamma_2$ at $\blacksquare$ (on the left-hand side of Figure \ref{fig:proof2_fig2}). Therefore, by Theorem \ref{thm:formula cv}, $F_{\mathcal{G}_{[\gamma]}}=F_{[\gamma]}$.
\end{itemize}
\item [2)] The case where $\tau_{[n-1]}$ is in the minimal perfect matching of $\hat{\mathcal{G}}_{\gamma_1}$ is analogous, exchanging the roles of $\gamma_1$ and $\gamma_2$. 

\end{itemize}

Finally, we prove that $\mathbf{g}_{\mathcal{G}_{[\gamma]}}=\mathbf{g}_{[\gamma]}$.  If $\operatorname{Res}([\gamma])=\{ \gamma_1 \}$, by construction, an edge with label $n$ is in $P_-(\mathcal{G}_{\gamma_1})$ if and only if two edges with label $n$ are in $P_-(\hat{\mathcal{G}}_{\gamma_1})$. Therefore, 
 
\begin{itemize}
    \item if the arc $\gamma_1$ does not cross $\tau_n$, then $\mathbf{g}_{\mathcal{G}_{[\gamma]}}=\mathbf{g}_{\hat{\mathcal{G}}_{\gamma_1}}=D\mathbf{g}_{\mathcal{G}_{\gamma_1}}=D\mathbf{g}_{\gamma_1}$;
    \item otherwise, if $\gamma_1$ crosses $\tau_n$, then $\mathbf{g}_{\mathcal{G}_{[\gamma]}}=\mathbf{g}_{\hat{\mathcal{G}}_{\gamma_1}}=D\mathbf{g}_{\mathcal{G}_{\gamma_1}}+\mathbf{e}_n=D\mathbf{g}_{\gamma_1}+\mathbf{e}_n$.
\end{itemize} 
On the other hand, if $\operatorname{Res}([\gamma])=\{ \gamma_1, \gamma_2 \}$, since $P_-(\mathcal{G}_{[\gamma]})$ is defined as the union of the minimal matchings of $\hat{\mathcal{G}}_{\gamma_1}$ and $\hat{\mathcal{G}}_{\gamma_2}$ (see Definition \ref{def:min matching}), and both $\gamma_1$ and $\gamma_2$ cross $\tau_n$, we have that
\begin{align*}
\mathbf{g}_{\mathcal{G}_{[\gamma]}}=\mathbf{g}_{\hat{\mathcal{G}}_{\gamma_1}}+\mathbf{g}_{\hat{\mathcal{G}}_{\gamma_2}}=D\mathbf{g}_{\mathcal{G}_{\gamma_1}}+D\mathbf{g}_{\mathcal{G}_{\gamma_2}}+2\mathbf{e}_n=D(\mathbf{g}_{\gamma_1}+\mathbf{g}_{\gamma_2}+\mathbf{e}_n).
\end{align*} 
In all cases, it follows from Theorem \ref{thm:formula cv} that  $\mathbf{g}_{\mathcal{G}_{[\gamma]}}=\mathbf{g}_{[\gamma]}$.   
\end{proof}

\section{Symmetric algebras from surfaces with a $\bbZ_2$-action}\label{s3}

In this section, given a skew-symmetrizable cluster algebra $\mathcal{A}_\bullet^\varsigma(\Tilde{T})$ with principal coefficients in an admissible $\varsigma$-invariant triangulation $\Tilde{T}$ of $(\Tilde{\textbf{S}},\Tilde{\textbf{M}})$, we associate a symmetric quiver algebra $A$ with it, in such a way that the non-initial cluster variables of $\mathcal{A}_\bullet^\varsigma(\Tilde{T})$ correspond to some orthogonal indecomposable $A$-modules. Building on the results of Section \ref{section:ssca_surfaces}, we then define a Caldero-Chapoton like map (see \cite{CC}) from these modules to the cluster algebra $\mathcal{A}_\bullet^\varsigma(\Tilde{T})$.

\vspace{0.3cm}

In the following, $k$ is an algebrically closed field. 
\subsection{Symmetric algebras and symmetric modules}\label{symmetric_quivers}
We first recall some essential notions of symmetric representation theory, introduced by Derksen and Weyman in \cite{DW}, as well as by Boos and Cerulli Irelli in \cite{BCI}, to set the notation.

\begin{definition}[Quiver and quiver algebra]
    For a quiver $Q=(Q_0,Q_1)$, where $Q_0$ is the finite set of vertices and $Q_1$ is the finite set of arrows, the \emph{path algebra} $kQ$ is the $k$-vector space generated by the set of all paths in $Q$, with multiplication given by concatenation of paths. Let $R$ be the two-sided ideal generated by the arrows of $Q$. Let $I \subseteq kQ$ be an \emph{admissible} ideal, that is, there exists an integer $m\geq 2$ such that $R^m \subseteq I \subseteq R^2$. Then the finite-dimensional quotient algebra $A=kQ/I$ is called a \emph{quiver algebra}. 
\end{definition}

In the following, any module is a right module. 

\begin{definition}[Symmetric quiver]\label{def_symm_quiver}
    A \emph{symmetric quiver} is a pair $(Q,\varrho)$, where $Q$ is a finite quiver and $\varrho$ is an involution of $Q_0$ and of $Q_1$ which reverses the orientation of arrows. 
\end{definition}

\begin{example}
Consider the following orientations of a Dynkin diagram of type $A_3$:
\begin{itemize}
    \item The quiver $Q = 1 \xrightarrow[]{a} 2 \xrightarrow[]{b} 3$ is symmetric, with $\varrho$ given by $\varrho(1)=3$, $\varrho(2)=2$ and $\varrho(a)=b$.
    \item The quiver $Q'= 1 \xrightarrow[]{a} 2 \xleftarrow[]{b} 3$ is not symmetric.
\end{itemize}
\end{example}

\begin{definition}[Symmetric quiver algebra]
    Let $(Q,\varrho)$ be a symmetric quiver. Let $I \subset kQ$ be an admissible ideal such that $\varrho(I)=I$. Then $A=kQ/I$ is called a \emph{symmetric quiver algebra}.
\end{definition}

\begin{definition}[Symmetric module]
    A \emph{symmetric module} over a symmetric algebra $A=kQ/I$ is a triple 
    $(V_i,\phi_a, \langle \cdot, \cdot \rangle)$, where $(V_i,\phi_a)$ is an $A$-module, $\langle \cdot, \cdot \rangle$ is a non-degenerate symmetric or skew-symmetric scalar product on $V=\displaystyle\bigoplus_{i \in Q_0}V_i$ such that its restriction to $V_i \times V_j$ is 0 if $j \neq \varrho(i)$, and $\langle \phi_a(v), w \rangle + \langle v, \phi_{\varrho(a)}(w) \rangle=0$, for every $a : i \to j \in Q_1$, $v \in V_i$, $w \in V_{\varrho(j)}$. If $\langle \cdot, \cdot \rangle$ is symmetric (resp. skew-symmetric), $(V_i,\phi_a, \langle \cdot, \cdot \rangle)$ is called $orthogonal$ (respectively, $symplectic$).
\end{definition}

\begin{definition}[Twisted dual]
    Let $L=(V_i,\phi_a)$ be a module over a symmetric algebra $A=kQ/I$. The \emph{twisted dual of $L$} is the $A$-module $\nabla L = (\nabla V_i, \nabla \phi_a)$, where $\nabla V_i=V_{\varrho(i)}^\ast$ and $\nabla \phi_a = - \phi_{\varrho (a)}^\ast$ ($\ast$ denotes the linear dual).
\end{definition}

\begin{remark}
The minus sign in the definition of the twisted dual is chosen so that the compatibility of the module structure with the scalar product $\langle\cdot,\cdot\rangle$ is equivalent to the existence of an isomorphism $L\cong\nabla L$. In particular, a module $L$ is symmetric if and only if $\langle \cdot, \cdot \rangle$ induces such an isomorphism. For further details, see \cite[Section 2.3]{BCI}.
\end{remark}


The following result, proved in \cite{DW} for symmetric quivers without relations, and in \cite{BCI} for any symmetric quiver algebra, shows that every indecomposable symmetric module is uniquely determined by the $\nabla$-orbit of an ordinary indecomposable module:

\begin{lemma}[{\cite[Proposition 2.7]{DW},\cite[Lemma 2.10]{BCI}}]\label{ind_symm}
 Let $N$ be an indecomposable symmetric module over a symmetric quiver algebra $A$. Then, one and only one of the following three cases can occur: 
 \begin{itemize}
     \item [(I)] $N$ is indecomposable as a $A$-module; in this case, $N$ is called of type (I), for “indecomposable”.
     \item [(S)] There exists an indecomposable $A$-module $L$ such that $N=L\oplus \nabla L$ and $L\ncong \nabla L$; in this case, $N$ is called of type (S), for “split”.
     \item [(R)] There exists an indecomposable $A$-module $L$ such that $N=L\oplus \nabla L$ and $L\cong \nabla L$; in this case, $N$ is called of type (R) for “ramified”.
 \end{itemize}
\end{lemma}

\subsection{From $\varsigma$-orbits to $\varrho$-orbits}

\begin{definition}[Quiver algebra associated with a triangulation]\label{def:quiver of triang}
    Let $T=\{\tau_1,\dots,\tau_n\}$ be a triangulation of a surface $(\mathbf{S},\mathbf{M})$. Let $B(T)=(b_{ij})$ be the signed adjacency matrix of $T$ (see Definition \ref{def:adj_matrix}). The \emph{quiver algebra associated with $T$} is the quiver algebra $A(T)=kQ(T)/I(T)$, where 
    \begin{itemize}
        \item $Q(T)$ is the quiver with vertices $1,\dots,n$, and for any $i \neq j$, $b_{ij}$ arrows from $j$ to $i$ if $b_{ij}>0$;
        \item $I(T)$ is generated by all paths $i\to j \to k$ such that there exists an arrow $k\to i$.
    \end{itemize}  
\end{definition}

\begin{remark}
    The algebra $A(T)$ is the Jacobian algebra of the quiver with potential associated in \cite{labardini2009quivers} with triangulations of unpunctured surfaces. Moreover, as shown in \cite{assem2010gentle}, $A(T)$ is a gentle algebra, and the arcs of $(\mathbf{S},\mathbf{M})$ correspond to strings and powers of bands of $A(T)$.
\end{remark}

Let $\Tilde{T}$ be an admissible $\varsigma$-invariant triangulation of a surface $(\Tilde{\textbf{S}},\Tilde{\textbf{M}})$ endowed with an orientation-preserving diffeomorphism $\varsigma$ of order 2. Then the algebra $A(\Tilde{T})$ associated with $\Tilde{T}$ is not symmetric.

\begin{example}
    Let $\Tilde{T}$ be the $\varsigma$-invariant triangulation in Figure \ref{fig:ex_adm_triang}. The algebra $A(\Tilde{T})$ is given by the quiver
    \begin{center}
\begin{tikzcd}
1 & 2 \arrow[l, "a"'] \arrow[d, "c"]   & 4 \arrow[l, "b"'] & 5 \arrow[l, "f"'] \arrow[r, "f'"] & 4' \arrow[r, "b'"] & 2' \arrow[r, "a'"] \arrow[d, "c'"]    & 1' \\
  & 3 \arrow[lu, "d"] \arrow[ru, "e"'] &                   &                                   &                    & 3' \arrow[lu, "e'"] \arrow[ru, "d'"'] &    
\end{tikzcd}
    \end{center}
    with relations $bc=ce=eb=0$ and $b'c'=c'e'=e'b'=0$. Clearly, there is only one non-trivial involution, which sends $i \mapsto i'$, for $i=1,\dots,4$ and fixes 5, and it does not reverse the orientation of the arrows. Thus, $A(\Tilde{T})$ is not symmetric.
\end{example}

\begin{figure}[h]
    \centering
    \includegraphics[width=0.9\linewidth]{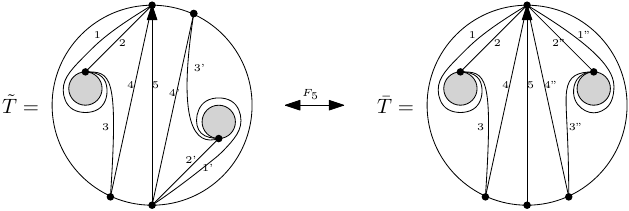}
    \caption{The operation $F_{\tau_n}$ applied to a planar surface with three boundary components.}
    \label{fig:f_d}
\end{figure}

\begin{figure}[h]
    \centering
    \includegraphics[width=1\linewidth]{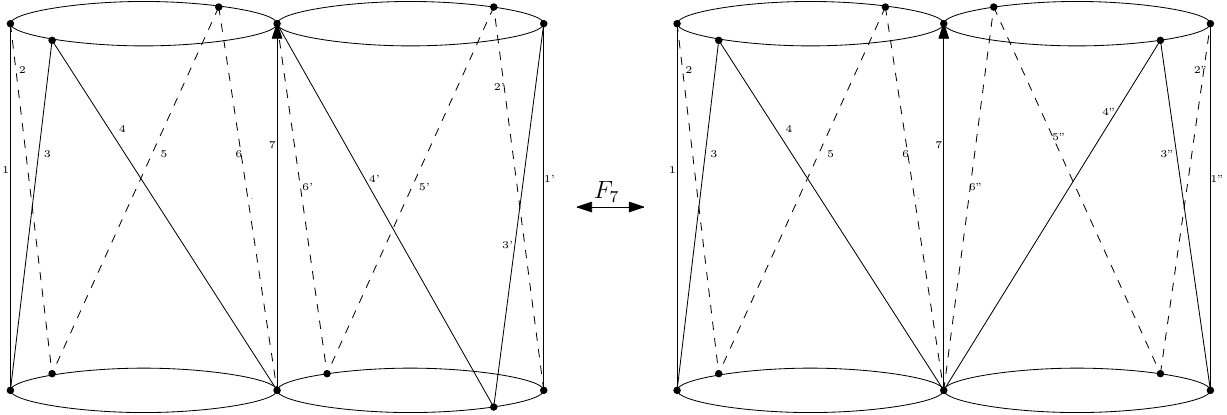}
    \caption{The operation $F_{\tau_n}$ applied to a non-planar surface formed by two cylinders glued along a vertical segment.}
    \label{fig:f_d_cilinders}
\end{figure}

In order to associate a symmetric algebra with a given admissible $\varsigma$-invariant triangulation $\Tilde{T}$ of $(\Tilde{\textbf{S}},\Tilde{\textbf{M}})$ with $\varsigma$-invariant arc $\tau_n$, we define an involution $F_{\tau_n}$ on the surface, depending on the orientation of $\tau_n$, as follows: 
\begin{definition}\label{def_f_d}
    $F_{\tau_n}$ is the operation on $(\Tilde{\textbf{S}},\Tilde{\textbf{M}})$ defined by the following two steps:
    \begin{itemize}
        \item [(1)] Cut the surface along $\tau_n$;
        \item [(2)] Replace the right part with a reflection of the left part along the axis of symmetry of $\tau_n$.
    \end{itemize}
    See Figures \ref{fig:f_d} and \ref{fig:f_d_cilinders} for two examples.
\end{definition}

\begin{figure}[h]
    \centering
    \includegraphics[width=0.8\linewidth]{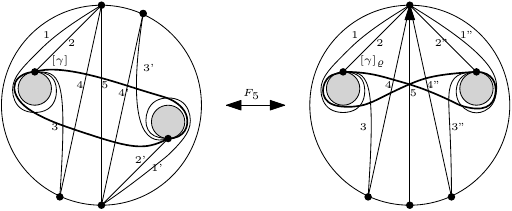}
    \caption{On the left, a $\varsigma$-orbit $[\gamma]$ in a triangulated surface of genus 0 with three boundary components; on the right, the $\varrho$-orbit $[\gamma]_\varrho$ corresponding to $[\gamma]$ in the flipped surface.}
    \label{fig:rho_orbit}
\end{figure}

\begin{remark}\label{remk:rho orbits}
The resulting surface $(\Bar{\textbf{S}},\Bar{\textbf{M}}):=F_{\tau_n}((\Tilde{\textbf{S}},\Tilde{\textbf{M}}))$ is endowed with an orientation-reversing diffeomorphism $\varrho$ of order 2. Furthermore, $F_{\tau_n}$ induces an action on the isotopy classes of the arcs of the polygon, in such a way that 
\[
\begin{array}{ccc}
\big\{\text{$\varsigma$-orbits of arcs of } (\Tilde{\textbf{S}},\Tilde{\textbf{M}})\big\} & \xlongleftrightarrow{F_{\tau_n}} & \big\{\text{$\varrho$-orbits of arcs of } (\Bar{\textbf{S}},\Bar{\textbf{M}})\big\}. \\
\text{$[\gamma]$} & & [\gamma]_\varrho
\end{array}
\] 
In particular, $\varsigma$-invariant arcs correspond to $\varrho$-invariant arcs; while $\varsigma$-invariant pairs of non-$\varsigma$-invariant arcs correspond to $\varrho$-invariant pairs of non-$\varrho$-invariant arcs. See Figure \ref{fig:rho_orbit} for an example.
\end{remark}

\begin{remark}\label{rmk_no_fixed_arrow}
Let $\Bar{T}:=F_{\tau_n}(\Tilde{T})=\{\tau_i\}_{i=1}^{n-1} \sqcup \{\tau_n\} \sqcup \{\tau_i''\}_{i=1}^{n-1}$. Then $\Bar{T}$ is an admissible $\varrho$-invariant triangulation of $(\Bar{\textbf{S}},\Bar{\textbf{M}})$. It follows that $A(\Bar{T})=kQ(\Bar{T})/I(\Bar{T})$ is a symmetric gentle algebra with respect to the involution induced by $\varrho$. Moreover, $Q(\Bar{T})$ has exactly one fixed vertex $n$, corresponding to $\tau_n$, and no fixed arrows. Since $A(\Bar{T})$ is a gentle algebra, indecomposable $A(\Bar{T})$-modules $L_\gamma$ correspond to arcs $\gamma$ of $(\Bar{\textbf{S}},\Bar{\textbf{M}})$ (see \cite{assem2010gentle} for more details). Furthermore, since $Q(\Bar{T})$ does not have fixed arrows, orthogonal indecomposable $A(\Bar{T})$-modules are either of type I or type S. It follows from Lemma \ref{ind_symm} that
\[
\begin{array}{ccc}
\big\{\varrho\text{-orbits of arcs of } (\Bar{\textbf{S}},\Bar{\textbf{M}})\text{ not in } \Bar{T}\big\} & \xlongleftrightarrow{} & \big\{\text{orthogonal indecomposable $A(\Bar{T})$-modules} \big\}. \\
\text{$[\gamma]_\varrho$} & & \text{$L_{[\gamma]_\varrho}$}
\end{array}
\] 
In particular, $\varrho$-invariant arcs $\gamma$ correspond to orthogonal indecomposable $A(\Bar{T})$-modules $L_{[\gamma]_\varrho}=L_\gamma$ of type I; while $\varrho$-invariant pairs $[\gamma]_\varrho=\{\gamma,\varrho(\gamma)\}$ of non-$\varrho$-invariant arcs correspond to orthogonal indecomposable $A(\Bar{T})$-modules $L_{[\gamma]_\varrho}=L_\gamma \oplus L_{\varrho(\gamma)}$ of type S.
\end{remark}

\begin{example}\label{ex:symm quiv triang}
 Let $\Bar{T}$ be the triangulation on the right-hand side of Figure \ref{fig:f_d}. The algebra $A(\Bar{T})$ is given by the quiver
    \begin{center}
        \begin{tikzcd}
1 & 2 \arrow[l, "a"'] \arrow[d, "c"]   & 4 \arrow[l, "b"'] & 5 \arrow[l, "f"'] & 4'' \arrow[l, "f''"'] \arrow[rd, "e''"'] & 2'' \arrow[l, "b''"'] & 1'' \arrow[l, "a''"'] \arrow[ld, "d''"] \\
  & 3 \arrow[lu, "d"] \arrow[ru, "e"'] &                   &                   &                                                  & 3'' \arrow[u, "c''"'] &                                                
\end{tikzcd}
    \end{center}
   with relations $bc=ce=eb=0$ and $b''e''=e''c''=c''b''=0$. Evidently, $A(\Bar{T})$ is symmetric with respect to the involution defined by $i \mapsto i''$, for $i=1,\dots,4$, and fixing 5.
\end{example}

\begin{definition}
    Let $L$ be an orthogonal indecomposable $A(\Bar{T})$-module. Then $L$ is called \emph{non-intersecting} (respectively, \emph{nice}) if $L=L_{[\gamma]_\varrho}$, and $F_{\tau_n}([\gamma]_\varrho)$ is a non-intersecting (respectively, nice) $\varsigma$-orbit.
\end{definition}

\begin{proposition}\label{prop:symm ind}
    Let $\Tilde{T}$ be an admissible $\varsigma$-invariant triangulation of a surface $(\Tilde{\textbf{S}},\Tilde{\textbf{M}})$ endowed with an orientation-preserving diffeomorphism $\varsigma$ of order 2, fixing globally $\Tilde{\textbf{M}}$. Let $\Bar{T}:=F_{\tau_n}(\Tilde{T})$. Then, we have the following bijection:
\[
\begin{array}{ccc}
\big\{\text{non-initial cluster variables of } \mathcal{A}_\bullet^\varsigma(\Tilde{T})\big\} & \xlongleftrightarrow{} & \big\{\text{non-intersecting $A(\Bar{T})$-modules} \big\}. \\
\text{$x_N$} & & \text{$N$}
\end{array}
\] 
\end{proposition}

\begin{proof}
     The result is obtained by combining Proposition \ref{prop:cluster var in A^sigma}, Remark \ref{remk:rho orbits} and Remark \ref{rmk_no_fixed_arrow}.
\end{proof}

In Section \ref{sect:symm ind to cl var}, we will see a purely representation-theoretic formula to express $x_N$ in terms of the initial cluster variables. 

\subsection{Cluster characters for gentle algebras}

In this section, we recall the definition of cluster character and the multiplication formula for gentle algebras established in \cite{ciliberti3}, which serves as a basis for Section \ref{sect:symm ind to cl var}.

\begin{definition}\label{def:g_vect}
    Let $A=kQ/I$ be a finite-dimensional algebra, and $L$ be an $A$-module. Let
    \begin{center}
        $0 \to L \xrightarrow[]{} I_0 \xrightarrow[]{} I_1$
    \end{center}
    be a minimal injective presentation of $L$, with $I_0=\displaystyle\bigoplus_{i\in Q_0}I(i)^{a_i}$ and $I_1=\displaystyle\bigoplus_{i\in Q_0}I(i)^{b_i}$, where $I(i)$ denotes the injective $A$-module at vertex $i$. Then the $\mathbf{g}$-vector of $L$ is the vector $\mathbf{g}_L \in \mathbb{Z}^{Q_0}$ whose $i$-th coordinate is given by
    \begin{center}
        $(\mathbf{g}_L)_i:=b_i-a_i$.
    \end{center}
\end{definition}  
\begin{definition}\label{defi::grassmannian}
Let $A$ be a finite-dimensional algebra. Let $\mathbf{e} \in \mathbb{Z}_{\geq 0}^{Q_0}$ and $L$ be an $A$-module. The \emph{quiver Grassmannian of $L$ with dimension vector $\mathbf{e}$} is the projective variety $\mathrm{Gr}_{\mathbf{e}}(L)$ of all submodules of $L$ of dimension vector $\mathbf{e}$.
\end{definition}
\begin{definition}\label{def:cc map}
        Let $A=kQ/I$ be a finite-dimensional algebra. Let $n=|Q_0|$, and let $B=B(Q)=(b_{ij})$ be the $n\times n$ matrix such that $b_{ij}=|\{a:j \to i \in Q_1\}|-|a:i\to j \in Q_1|$. Let $L$ be an $A$-module. The \emph{cluster character of $L$}, also known as \emph{Caldero-Chapoton map}, is the Laurent polynomial
    \begin{center}
        $CC(L)=\displaystyle\sum_{\mathbf{e}\in \mathbb{Z}_{\geq 0}^n}\chi(\mathrm{Gr}_{\mathbf{e}}(L))\mathbf{x}^{\,B\mathbf{e}+\mathbf{g}_L}\,\mathbf{y}^{\mathbf{e}}\in \mathbb{Z}[y_1,\dots,y_n][x_1^{\pm 1},\dots,x_n^{\pm 1}]$,
    \end{center}
    where $\chi$ is the Euler-Poincar\'e characteristic.
    The $F$-polynomial of $L$ is defined to be the specialization of $CC(L)$ at $x_1=\cdots=x_n=1$. In other words,
     $$ F_L(\mathbf{y}) := \displaystyle\sum_{\mathbf{e}\in\mathbb{Z}_{\geq 0}^{Q_0}} \chi\big( \mathrm{Gr}_{\mathbf{e}}(L) \big) \mathbf{y}^{\mathbf{e}}.
 $$
\end{definition}

\begin{definition}
Let $A$ be a finite-dimensional algebra, and let $M,N$ be $A$-modules. We say that $M \leq_{\mathrm{Ext}} N$ if there exist $A$-modules $M_1,\dots, M_k$ such that for every $i$ there exists a short exact sequence
\begin{center}
    $0 \to U_i \to M_{i-1} \to V_i \to 0$,
\end{center}
such that $M_1=M$, $M_k=N$, $M_i \cong U_i \oplus V_i$.
\end{definition}

\begin{theorem}[{\cite[Theorem 4.0.11]{ciliberti3}}] \label{the:cc mult formula}
    Let $A(T)$ be the gentle algebra associated with a triangulation $T$ of a surface $(S,M)$. Let $X$, $S$ be rigid (that is, $\operatorname{dim} \operatorname{Ext}^1(X,X)=\operatorname{dim} \operatorname{Ext}^1(S,S)=0$) and indecomposable $A$-modules such that $\operatorname{dim} \operatorname{Ext}^1(S,X)=1$. Let $\xi \in \mathrm{Ext}^1(S,X)$ be a non-split short exact sequence with middle term $Y$. Then, denoting by $\overline{X}$ the kernel of a non-zero morphism from $X$ to $\tau S$ that does not factor through an injective $A$-module, and by $\underline{S}$ the image of a non-zero morphism from $\tau^{-1} X$ to $S$ that does not factor through a projective $A$-module,
\begin{equation}\label{eq_cor1}
        CC(X)CC(S)= CC(Y)\bold{x}^{\bold{g}_X+\bold{g}_S-\bold{g}_Y} + \bold{y}^{\textbf{dim} \underline{S}}CC(M)\bold{x}^{B\textbf{dim}\underline{S}+\bold{g}_X+\bold{g}_S-\bold{g}_M},
\end{equation}
where $M$ is the $\leq_{\mathrm{Ext}}$-minimum extension between $S/\underline{S}$ and $\overline{X}$. Moreover, \eqref{eq_cor1} is an exchange relation between the cluster variables $CC(X)$ and $CC(S)$ in the cluster algebra $\mathcal{A}_\bullet(T)$ with principal coefficients in $T$.   
\end{theorem}

\begin{remark}
    The existence and uniqueness of the $\leq_{\mathrm{Ext}}$-minimum extension between $S/\underline{S}$ and $\overline{X}$ is part of \cite[Theorem 4.0.11]{ciliberti3}.
\end{remark}

\begin{remark}\label{rmk_spec}
    Specializing at $x_1=\cdots=x_n=1$, where $n=|Q_0|$, we get the following multiplication formula for $F$-polynomials:
     \begin{equation}\label{eq_rmk_cor1}
        F_XF_S= F_Y + \bold{y}^{\textbf{dim} \underline{S}}F_M.
    \end{equation} 
\end{remark}

\subsection{Cluster variables corresponding to orthogonal indecomposable modules}\label{sect:symm ind to cl var}
In this section, we work in the following setting:
\begin{itemize}
    \item $\Tilde{T}$ is an admissible $\varsigma$-invariant triangulation of a surface $(\Tilde{\mathbf{S}},\Tilde{\mathbf{M}})$ that admits an orientation-preserving diffeomorphism $\varsigma$ of order 2;
    \item $\Bar{T}=F_{\tau_n}(\Tilde{T})$ is the triangulation of the flipped surface $(\Bar{\mathbf{S}},\Bar{\mathbf{M}})$ that admits an non-orientation-preserving diffeomorphism $\varrho$ of order 2;
    \item $T=\operatorname{Res}(\Bar{T})=\operatorname{Res}(\Tilde{T})=\{\tau_1, \dots, \tau_n\}$ is the triangulation of the collapsed surface $(\mathbf{S},\mathbf{M})$.
\end{itemize}
The restriction on non-intersecting $\varsigma$-orbits corresponds to the following operation on non-intersecting $A(T)$-modules: 
\begin{definition}\label{def_res}
\begin{itemize}
    \item [(i)] Let $N=(V_i,\phi_a,\langle \cdot, \cdot \rangle)$ be an non-intersecting $A(\Bar{T})$-module. Then the \emph{restriction} of $N$ is the $A(T)$-module $\operatorname{Res}(N)=(\operatorname{Res}(V)_i,\operatorname{Res}(\phi)_a)$, where $\operatorname{Res}(V)_i=V_i$ for any $i\in Q(T)_0$ and $\operatorname{Res}(\phi)_a=\phi_a$ for any $a \in Q(T)_1$.
    \item [(ii)] Let $v \in \mathbb{Z}_{\geq 0}^{2n-1}$. The \emph{restriction} of $v$, denoted by $\operatorname{Res}(v)$, is the vector $\operatorname{Res}(v) \in \mathbb{Z}_{\geq 0}^{n}$ of the first $n$ coordinates of $v$.
\end{itemize}
\end{definition}

For an non-intersecting $A(\Bar{T})$-module $N$, $F_N$ and $\mathbf{g}_N$ denote the $F$-polynomial and the $\mathbf{g}$-vector of the non-initial cluster variable $x_N$ of $\mathcal{A}_\bullet^\varsigma(\Tilde{T})$ that corresponds to $N$ by Proposition \ref{prop:symm ind}. On the other hand, $F_{\operatorname{Res}(N)}$ and $\mathbf{g}_{\operatorname{Res}(N)}$ are the $F$-polynomial and the $\mathbf{g}$-vector of the $A(T)$-module $\operatorname{Res}(N)$, as in Definitions \ref{def:cc map} and \ref{def:g_vect}. The following theorem gives us a purely representation-theoretic formula to compute $F_N$ and $\mathbf{g}_N$:
\begin{theorem}\label{cat_interpr}
Let $N$ be a nice $A(\Bar{T})$-module. Let $D=\operatorname{diag}(1,\dots,1,2)\in \mathbb{Z}^{n\times n}$.
\begin{itemize}
    \item [(i)] If $\operatorname{Res}(N)=(V_i,\phi_a)$ is indecomposable as $A(T)$-module, then 
    \begin{equation*}\label{e_1}
        \text{$F_N=F_{\operatorname{Res}(N)}$,}
    \end{equation*}
    and
    \begin{equation*}
 \mathbf{g}_{N}=\begin{cases}
          \text{$D \mathbf{g}_{\operatorname{Res}(N)}$ \hspace{1.8cm}if $\operatorname{dim} V_n =0$;}\\
          \text{$D \mathbf{g}_{\operatorname{Res}(N)}+\mathbf{e}_n$ \hspace{1cm}if $\operatorname{dim} V_n \neq 0$.}
      \end{cases}
      \end{equation*}
\item [(ii)]      Otherwise, $N=L\oplus \nabla L$ with  $\operatorname{dim} \operatorname{Ext}^1(\nabla L, L)=1$, and then there exists a non-split short exact sequence
    \begin{equation*}
        0 \to L \to G_1 \oplus G_2 \to \nabla L \to 0,
    \end{equation*}
    where $G_1$ and $G_2$ are orthogonal indecomposable $A(\Bar{T})$-modules of type I. Furthermore,
    \begin{equation*}
        F_N= F_{\operatorname{Res}(N)} - \mathbf{y}^{\operatorname{Res}(\textbf{dim}\underline{\nabla L})}F_{\operatorname{Res}(M)},
    \end{equation*}
    and
    \begin{equation*}
        \mathbf{g}_{N}=D(\mathbf{g}_{\operatorname{Res}(N)}+\mathbf{e}_n),
    \end{equation*}
where $M$ is the $\leq_{\mathrm{Ext}}$-minimum extension in $A(\Bar{T})$ between $\nabla L/\underline{\nabla L}$ and $\overline{L}$.
 \end{itemize}  
     
\end{theorem}

\begin{remark}
    If $(\Tilde{\textbf{S}},\Tilde{\textbf{M}})$ is a regular polygon, we recover \cite[Theorem 5.0.14]{ciliberti3}.
\end{remark}

\begin{proof}[Proof of Theorem \ref{cat_interpr}]
Let $[\gamma]_\varrho$ be the $\varrho$-orbit corresponding to $N$. If $\operatorname{Res}(N)$ is indecomposable as $A(T)$-module, then $\operatorname{Res}([\gamma])=\operatorname{Res}([\gamma]_\varrho)=\{\gamma_1\}$. Thus,
\begin{align*}
    F_{\operatorname{Res}(N)}=F_{\gamma_1}=F_{[\gamma]}=F_N,
\end{align*}
where the second last equality is given by Theorem \ref{thm:formula cv} (i).
Assume now that $\operatorname{Res}(N)$ is not indecomposable as a $A(T)$-module. It follows that $[\gamma]_\varrho=\{\gamma,\gamma''\}$, where $\gamma$ and $\gamma''$ are two arcs of $(\Bar{\textbf{S}},\Bar{\textbf{M}})$ that cross exactly once at a point $x \in \tau_n$. 
\begin{figure}[h]
    \centering
    \includegraphics[width=0.7\linewidth]{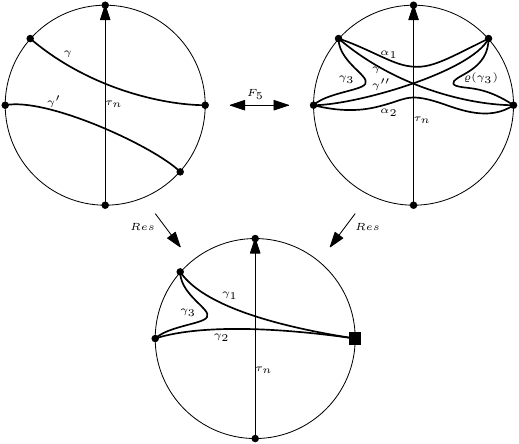}
    \caption{Illustration of the proof of Theorem \ref{cat_interpr} (ii).}
    \label{fig:proof_thm3}
\end{figure}
By \cite[Theorem 3.7]{CSA}, 
$\operatorname{dim} \operatorname{Ext}^1(L_{\gamma''}, L_{\gamma})=1$ and, in the notation of Figure \ref{fig:proof_thm3}, there is a non-split short exact sequence
\begin{align*}
    0 \to L_{\gamma} \to L_{\alpha_1}\oplus L_{\alpha_2} \to L_{\gamma''} \to 0.
\end{align*}

Since $\alpha_1$ and $\alpha_2$ are $\varrho$-invariant arcs by construction, $G_1:=L_{\alpha_1}$ and $G_2:=L_{\alpha_2}$ are orthogonal indecomposable $A(\Bar{T})$-modules of type I. Furthermore, by Theorem \ref{the:cc mult formula} and Remark \ref{rmk_spec} applied to $L:=L_{\gamma}$ and $\nabla L= L_{\gamma''}$, in the cluster algebra associated with $(\Bar{\textbf{S}},\Bar{\textbf{M}})$ with principal coefficients in $\Bar{T}$, the following relation holds:
\begin{align*}
    F_{L\oplus \nabla L}&=F_{G_1\oplus G_2}+\mathbf{y}^{\textbf{dim}\underline{\nabla L}}F_M,
\end{align*}
where $M$ is the extension between $\nabla L/\underline{\nabla L}$ and $\overline{L}$ in $A(\Bar{T})$ that is minimal with respect to the $\operatorname{Ext}$-order.
On the other hand, by Proposition \ref{up:skein1},
\begin{align*}
  F_{L\oplus \nabla L}&=F_{G_1\oplus G_2}+\mathbf{y}^{\mathbf{d}_{\alpha_1,\alpha_2}}F_{L_{\gamma_3}\oplus L_{\varrho(\gamma_3)}}.  
\end{align*}
Thus,
\begin{align*}
  \textbf{dim}\underline{\nabla L}&=\mathbf{d}_{\alpha_1,\alpha_2},   
\end{align*}
and
\begin{align*}
   M&=L_{\gamma_3}\oplus L_{\varrho(\gamma_3)}.
\end{align*}
If $\operatorname{Res}([\gamma])=\operatorname{Res}([\gamma]_\varrho)=\{\gamma_1,\gamma_2\}$, then 
\begin{align*}
    F_{\operatorname{Res}(N)} - \mathbf{y}^{\operatorname{Res}(\textbf{dim}\underline{\nabla L})}F_{\operatorname{Res}(M)}&=F_{L_{\gamma_1}}F_{L_{\gamma_2}}-\mathbf{y}^{\mathbf{d}_{\gamma_1,\gamma_2}}F_{L_{\gamma_3}}\\&=F_{\gamma_1}F_{\gamma_2}-\mathbf{y}^{\mathbf{d}_{\gamma_1,\gamma_2}}F_{\gamma_3}\\&=F_{[\gamma]}=F_N,
\end{align*}
where the second last equality is given by Theorem \ref{thm:formula cv} (ii). Similarly,
\begin{align*}
D(\mathbf{g}_{\operatorname{Res}(N)}+\mathbf{e}_n)&=D(\mathbf{g}_{L_{\gamma_1}}+\mathbf{g}_{L_{\gamma_2}}+\mathbf{e}_n)\\&=D(\mathbf{g}_{\gamma_1}+\mathbf{g}_{\gamma_2}+\mathbf{e}_n)\\&=\mathbf{g}_{[\gamma]}=\mathbf{g}_{N}.
\end{align*}
\end{proof}

\begin{example}
    Let $\mathcal{A}_\bullet^\varsigma(\Tilde{T})$ be the skew-symmetrizable cluster algebra with principal coefficients in the triangulation $\Tilde{T}$ in Figure \ref{fig:ex_adm_triang}. Let $A(\Bar{T})$ be the corresponding symmetric algebra described in Example \ref{ex:symm quiv triang}. We consider the nice $A(\Bar{T})$-module $N=\begin{smallmatrix}\hspace{0.1cm}3\hspace{0.05cm}5\\\hspace{-0.1cm}1\hspace{0.05cm}4 \end{smallmatrix}\oplus \begin{smallmatrix}\hspace{0.1cm}4''\hspace{0.05cm}1''\\\hspace{-0.15cm}5\hspace{0.07cm}3'' \end{smallmatrix}$. Let $x_N$ be the cluster variable of $\mathcal{A}_\bullet^\varsigma(\Tilde{T})$ that corresponds to $N$, and let $F_N$ and $\mathbf{g}_N$ denote its $F$-polynomial and its $\mathbf{g}$-vector, respectively. By Theorem \ref{cat_interpr},
    \begin{align*}
        F_N&=F_{\operatorname{Res}(N)}-y_5F_{\operatorname{Res}(\begin{smallmatrix}
            1
        \end{smallmatrix}\oplus \begin{smallmatrix}
            1''
        \end{smallmatrix})}=F_{\begin{smallmatrix}\hspace{0.1cm}3\hspace{0.05cm}5\\\hspace{-0.1cm}1\hspace{0.05cm}4 \end{smallmatrix}\oplus \begin{smallmatrix}
            5
        \end{smallmatrix}}-y_5F_{\begin{smallmatrix}
            1
        \end{smallmatrix}}\\&=y_1y_3y_4y_5^2 + 2y_1y_3y_4y_5 + y_1y_4y_5^2 + y_1y_3y_4 + 2y_1y_4y_5 + y_4y_5^2 + y_1y_4 + 2y_4y_5 + y_1 + y_4 + 1,\\
        \mathbf{g}_N&= D(\mathbf{g}_{\operatorname{Res}(N)}+\mathbf{e}_n)=D(\mathbf{g}_{\begin{smallmatrix}\hspace{0.1cm}3\hspace{0.05cm}5\\\hspace{-0.1cm}1\hspace{0.05cm}4 \end{smallmatrix}\oplus \begin{smallmatrix}
            5
        \end{smallmatrix}} + \mathbf{e}_5)=\begin{pmatrix}
            -1\\1\\1\\-1\\0
        \end{pmatrix}.
    \end{align*}
 
\end{example}

\section{Acknowledgements}
I would like to thank the anonymous referee for their valuable comments and suggestions, which have significantly improved the presentation of the paper.

\bibliographystyle{amsalpha-fi-arxlast}
\bibliography{reference.bib}

\end{document}